\magnification 1200
\input amssym.def
\input amssym.tex
\parindent = 40 pt
\parskip = 12 pt
\font \heading = cmbx10 at 12 true pt
 at 22 true pt
 at 19 true pt
 at 7 true pt
\def \R{{\bf R}}

\centerline{\heading Maximal averages over hypersurfaces and the Newton polyhedron}
\rm
\line{}
\line{}
\centerline{\heading Michael Greenblatt}
\line{}
\centerline{February 7, 2010}
\baselineskip = 12 pt
\font \heading = cmbx10 at 14 true pt
\line{}
\line{}
\noindent{\bf 1. Introduction and statement of results}

\vfootnote{}{This research was supported in part by NSF grant DMS-0919713} In this paper we are concerned
with two closely related objects in analysis, maximal averages over hypersurfaces and Fourier transforms of
surface-supported measures. Let $S$ be a smooth hypersurface in $\R^{n+1}$ for 
$n \geq 2$, and let $\sigma$ denote the standard surface measure on $S$.  
We consider the maximal operator, defined initially on Schwarz functions, given by
$$Mf(x) = \sup_{t > 0} |\int_S f(x - ts) \phi(s)\,d\sigma(s)| \eqno (1.1)$$
Here $\phi(x)$ is a smooth cutoff function that localizes the surface $S$ near some specific $y \in S$.
The goal here is to determine the values of $p$ for which $M$ is bounded on $L^p$. The earliest work on
this subject was done in the case where $S$ is a sphere, when Stein [St1] showed $M$ is bounded on $L^p$ 
iff $p > {n + 1 \over n}$ for $n > 1$. This was later generalized by Greenleaf [Gr] to surfaces of 
nonvanishing Gaussian curvature, and the $n = 1$ case was later proven by Bourgain [Bo]. Since
then, there have been a wide range of papers on this subject, which we will describe in more detail 
throughout this section. Although there are
many interesting issues when $p \leq 2$, for the purposes of this paper we always assume $p > 2$. Note that
if $M$ is bounded on some $L^p$, by interpolating with the $L^{\infty}$ case one has that $M$ is bounded on
$L^{p'}$ for $p' > p$. Hence our goal is to determine the optimal $p_0 \geq 2$ for which $M$ is bounded on 
$L^p$ for $p > p_0$. 

Let $L$ be an invertible linear transformation, and let $M_L$ be the maximal operator corresponding to the
surface $L(S)$. Then one can easily check from the definitions that $M_L f (x) = |det(L)| M(f \circ L)
(L^{-1}x)$. Hence in our arguments we may replace $M$ by $M_L$ at will. In particular, without loss of
generality we henceforth assume that $(0,...,0,1)$ is not in the tangent plane $T_y(S)$. We do this so that 
we may represent $S$ near $y$ as the graph of some function $g(x_1,...,x_n)$, which permits us to do the 
coordinate-dependent analysis of this paper. Let $z$ be such that $(z,g(z)) = y$ and define $h(x) = 
g(x) - g(z) - \nabla g(z) \cdot (x - z)$. Geometrically, $h(x)$ is the vertical distance from 
$S$ to $T_y(S)$ over $x$. If $y \notin T_y(S)$
and $h(x)$ has a zero of infinite order at $z$, it is not hard to show that $M$ is bounded on no $L^p$ 
space for $p < \infty$ as long as $\phi$ is nonnegative with $\phi(y) \neq 0$. On the other hand, by a 
result of Sogge-Stein [SoSt] if the Gaussian curvature of $S$ does not 
vanish to infinite order at $y$ then the reverse holds; $M$ is bounded on $L^p$ for some finite $p$. (See
also [CoMa] for another theorem of this kind). In Corollary 1.4 we will generalize this further.

The other (related) subject we are interested in this paper is the decay of Fourier transforms of surface 
measures. Let $S$ and $\phi$ be as above and consider $T(\lambda)$ defined by
$$T(\lambda) = \int_S e^{-i \lambda \cdot x} \phi(x)\,d\sigma(x) \eqno (1.2)$$
$T$ may be recognized as the Fourier transform of the surface measure of $S$ localized around $y$. As is
well known, $T$ is closely related to the maximal operator $M$. We may assume that the support of $\phi(x)$ 
is small enough such that we may write
$$T(\lambda) = \int e^{-i\lambda_1 x_1 - ... - i\lambda_n x_n - i\lambda_{n+1}g(x_1,...,x_n)}
\psi(x_1,...,x_n)\,dx_1...\,dx_n \eqno (1.3)$$
Here $\psi(x_1,...,x_n)$ is now a cutoff function localized around the $z$ such that $(z,g(z)) = y$. The main
goal for $T(\lambda)$ is to determine the optimal $\epsilon > 0$ for which one has an estimate
$$|T(\lambda)| \leq C|\lambda|^{-\epsilon} \eqno (1.4)$$
Clearly, by shifting coordinates, without loss of generality we may assume that $z = 0$ and replace $g(x)$
by $G(x) = g(x + z)$ and $\psi(x)$ by $\Psi(x) = \psi(x + z)$. Note that up to a 
magnitude one factor, $T(\lambda)$ is equal to
$$\int e^{-i\lambda_1x_1... - i\lambda_nx_n  - i\lambda_{n+1}(G(x_1,...,x_n) - G(0,...,0))}
\Psi(x_1,...,x_n)\,dx_1...\,dx_n \eqno (1.5)$$
Furthermore, one may write $G(x_1,...,x_n) - G(0,...,0)  =  \sum_{i=1}^n c_ix_i + H(x_1,...,x_n)$, where 
$c_i$ denotes $\partial_{x_i}G(0,...,0)$ and where $H(x)$ has a zero of order at least 2 at the origin. We
rewrite $(1.5)$ as
$$\int e^{-i(\lambda_1 + c_1 \lambda_{n+1})x_1... - i(\lambda_n + c_n \lambda_{n+1})x_n  - i\lambda_{n+1}
H(x_1,...,x_n)} \Psi(x_1,...,x_n)\,dx_1...\,dx_n \eqno (1.6)$$
Since invertible linear transformations in the $\lambda$ variables do not affect the truth of $(1.4)$, one 
may replace $\lambda_i$ for $i \leq n$ by $\lambda_i + c_i \lambda_{n+1}$ and assume $T(\lambda)$ is given 
by 
$$\int e^{-i\lambda_1 x_1... - i\lambda_n x_n  - i\lambda_{n+1} H(x_1,...,x_n)} \Psi(x_1,...,x_n)\,dx_1...\,
dx_n \eqno (1.7)$$
Stated another way, in the analysis of the Fourier transforms of surface-supported measures with regards to
$(1.4)$, one may always replace $g(x)$ by  
$H(x) = g(x+z) - g(z) - \nabla g (z) \cdot x$, a function with a critical point at the origin. Geometrically,
$H(x)$ is vertical distance from $T_y(S)$ to $S$ over the point $x + z$. Correspondingly,
our theorems for Fourier transforms of surface-supported measures will be stated in terms of $H(x)$. Because
of the close ties between those theorems and our theorems about maximal operators, our theorems on the
maximal operators will also be given in terms of $H(x)$. In both cases, we will be applying the resolution of
singularities methods of [G1] to $H(x)$. For the maximal operators, we will use them in
conjunction with appropriate damping methods analogous to those used in a number of papers in this subject 
such as
[SoSt]. We will obtain sharp estimates for a large class of maximal operators that includes for example
the convex surfaces of finite line type treated in several papers such as  [IoSa1] [IoSa2] [IoSaSe] and
[NaSeWa]. We will prove the analogous results for the Fourier transforms of surface measures, which will 
similarly include the convex surfaces of finite line type considered by [BrNaWa]. We now give some terminology 
used throughout this paper. 

\noindent {\bf Definition 1.1.} Let $f(x)$ be a smooth function defined on a neighborhood of the origin in 
$\R^n$, and
let $f(x) = \sum_{\alpha} f_{\alpha}x^{\alpha}$ denote the Taylor expansion of $f(x)$ at the origin.
For any $\alpha$ for which $f_{\alpha} \neq 0$, let $Q_{\alpha}$ be the octant $\{x \in \R^n: 
x_i \geq \alpha_i$ for all $i\}$. Then the {\it Newton polyhedron} $N(f)$ of $f(x)$ is defined to be 
the convex hull of all $Q_{\alpha}$.  

A Newton polyhedron can contain faces of various dimensions in various configurations. 
These faces can be either compact or unbounded. In this paper, as in earlier work like [G1] and [V], an 
important role is played by the following functions, defined for compact faces of
the Newton polyhedron. A vertex is always considered to be a compact face of dimension zero.

\noindent {\bf Definition 1.2.} Suppose $F$ is a compact face of $N(f)$. Then
if $f(x) = \sum_{\alpha} f_{\alpha}x^{\alpha}$ denotes the Taylor expansion of $f$ like above, 
define $f_F(x) = \sum_{\alpha \in F} f_{\alpha}x^{\alpha}$.

\noindent In order to state our theorems we use the following terminology.

\noindent {\bf Definition 1.3.} Assume $N(f)$ is nonempty. Then the 
{\it Newton distance} $d(f)$ of $f(x)$ is defined to be $\inf \{t: (t,t,...,t,t) \in N(f)\}$.

Our results will be sharp when the order of any zero of any $H_F(x)$ on $(\R - \{0\})^n$ is at most 
$d(H)$, where $H(x)$ is as in the previous section and $H_F(x)$ is as in Definition 1.2. This is also
where the sharpest results in [G1] occur. Roughly speaking, under these conditions, for the
purposes of this paper the zero of $H(x)$ at the origin is stronger than than any other zero of $H(x)$ in a 
neighborhood of the origin. A more precise way of saying this in the real-analytic case is that there is a small ball $B$ centered at the 
origin such that if $\int_B |H|^{-\delta}$ infinite, then for any sufficiently small 
neighborhood $U$ of the origin one also has $\int_U |H|^{-\delta}$ is infinite, regardless of how small $U$
is. 

Similar to [G1], we will also have partial results for when the maximum order of a zero of some $H_F(x)$ 
on $(\R - \{0\})^n$ is 
greater than $d(H)$. While these results are usually not sharp, they do provide
easily stated upper bounds for the minimal $p_0$ for which $M$ is bounded on $L^p$ for $p > p_0$ (we will 
see that $p_0$ can be taken to be the order of the zero of $H$ at the origin), and 
similarly for the supremum of the $\epsilon$ for which $(1.4)$ holds. They also will be used in the proof of
the extension of [SoSt] given by Corollary 1.4.

We now come to our theorems. In the following, we assume like above that we are in coordinates such that
$(0,...,0,1) \notin T_y(S)$ so that near $y$, the surface $S$ is the graph of a smooth function $g(x)$. 
We let $z$ be such that $(z,g(z)) = y$, and as indicated above, our theorems will be stated in terms of 
$H(x)= g(x+z) - g(z) - \nabla g(z) \cdot x$, a function which has a critical point at the origin. The 
Newton polyhedron $N(H)$ is useful in understanding the singularity of $H(x)$ at the origin, and our main 
theorem will be stated in terms of this Newton polyhedron. 

\noindent {\bf Theorem 1.1.} Suppose that $g(x)$ is real-analytic and that the Hessian matrix of $g(x)$ 
does not have rank $\leq 1$ everywhere. Then there 
is a neighborhood $U$ of $y$ such that if $\phi$ is supported in $U$ the following hold. 

\noindent {\bf a)} If $d(H) \leq 2$ and each zero of each $H_F(x)$ on $(\R - \{0\})^n$ has order at
most 2, then $M$ is bounded on $L^p$ for all $p > 2$. 

\noindent {\bf b)} If $d(H) > 2$ and each zero of each $H_F(x)$ on $(\R - \{0\})^n$ has order at most
$d(H)$, then $M$ is bounded on $L^p$ for all $p > d(H)$. If in addition $y \notin T_y(S)$ and $\phi$ is 
nonnegative with $\phi(y) \neq 0$, then $M$ is unbounded on $L^p$ for $p \leq d(H)$. 

\noindent {\bf c)} If the maximum order of any zero of any $H_F(x)$ on $(\R - \{0\})^n$ is given by
$m > \max(d(H),2)$, then $M$ is bounded on $L^p$ for all $p > m$.

Note that any Newton polyhedron falls into one of the parts of Theorem 1.1 (although the Hessian
condition may not be satisfied). In Lemma 4.1, we will see that in case b) of Theorem 1.1, the Hessian
condition can only be not satisfied if $N(H)$ has exactly one vertex, and that vertex lies on a coordinate
axis.

The Hessian condition  can be
understood as follows. It is well-known (and not hard to show directly) that if the Hessian of $g(x)$ does 
have rank $\leq 1$ everywhere and $g$ is a polynomial then the level sets of $H$ are hyperplanes.
Since these hyperplanes have to be parallel in order to not intersect, after
a linear coordinate change $H$ is a function of one variable. The situation now resembles that of a
maximal operator over curves in $\R^2$, and one might expect from [Bo] [Io] [MocSeSo] that the analysis
of such maximal operators may require additional ideas. For a general real-analytic function $g(x)$
the level sets aren't as easily described, but these situations will still fall out of the scope of this 
paper, perhaps for similar reasons. 

In part b) of the above theorem, we stipulate that $y \notin T_y(S)$ for
sharpness so that one can invoke a result of Iosevich and Sawyer in [IoSa1] giving a necessary condition for
boundedness of $M$ on $L^p$. It is unclear what conditions might be sufficient and necessary for $M$ to be bounded $L^p$
for all $p \leq m$ in the setting of Theorem 1.1c). However, the next result gives some situations in which
it can be readily shown that Theorem 1.1c) is sharp. 

\noindent {\bf Theorem 1.2.} Once again suppose $g(x)$ is real-analytic. Suppose $y \notin T_y(S)$ and
$\phi$ is nonnegative with $\phi(y) \neq 0$. Let $m$ denote the maximum order of any zero of any $H_F(x)$ on
$(\R - \{0\})^n$ and suppose $m > \max(d(H),2)$. Suppose there is a sequence $a^k$ of points in
$(\R - \{0\})^n$ with $a^k \rightarrow 0$ such that on a neighborhood of each $a^k = (a^k_1,...,a^k_n)$ 
the function $H(x)$ is of the form 
$\xi_k(x)(x_i - b_k(x_1,...x_{i-1},x_{i+1},...,x_n))^m$ for some real-analytic $b_k$ and $\xi_k$ such that
$\xi_k(a^k) \neq 0$ and $b_k((a^k_1,...,a^k_{i-1},a^k_{i+1},...,a^k_n)) = a^k_i$. Then $M$ is 
unbounded on $L^p$ for all $p \leq m$.

We use the theorem of [IoSa1] to prove Theorem 1.2 in short order. This theorem implies that if 
$\int |H|^{-{1 \over p}}$ is infinite on every neighborhood of the origin then $M$ is unbounded on $L^p$ as
long as $y \notin T_y(S)$. Clearly if $H(x)$ is of the form of Theorem 1.2, then the integral of $|H(x)|^
{-{1 \over m}}$ is infinite on a neighborhood of each $a^k$ and therefore on any neighborhood of the origin. 
Hence as long as $y \notin T_y(S)$, $M$ is unbounded on $L^m$ and the exponent $m$ of Theorem 1.1c) is 
sharp. 

The next theorem extends Theorem 1.1 to a class of smooth surfaces that will be seen to include the convex 
surfaces of finite line type considered by several authors.

\noindent {\bf Theorem 1.3.} Suppose $g(x)$ is a smooth function such that $N(H)$ intersects each  
coordinate axis. Then the conclusions of Theorem 1.1 hold as long as there are
directions $u$ and $v$ such that the Hessian determinant of $g(x)$ in the $u$ and $v$ variables, viewed
as a function of $x$, does not vanish to infinite order at $z$.  

\noindent {\bf Examples.}

\noindent {\bf Example 1.} Suppose $g$ is real-analytic and $N(H)$ is nondegenerate in the sense of Varchenko.
This means that any zero of any $H_F(x)$ on $(\R - \{0\})^n$ has order at most 1. By Lemma 4.1, the Hessian 
condition of Theorem 1.1 is automatically satisfied unless $N(H)$ has exactly one vertex and it lies on a
coordinate axis. As a result, other than in these exceptional cases Theorem 1.1a-b says that if the 
Newton distance $d(H)$ of $H$ is at most 2, $M$ is bounded on $L^p$ for all $p > 2$, and that
if $d(H) > 2, M$ is bounded on $L^p$ for $p$ greater than $d(H)$, with this $p$ optimal as long as 
$y \notin T_y(S)$. 

\noindent {\bf Example 2.} Suppose $S$ is a convex surface of finite line type as considered in [IoSa1] 
[IoSa2] [IoSaSe] [NaSeWa] and others. Then by a theorem of Schulz [Sc] there is a linear coordinate change 
after which the Newton polyhedron of $H$ has exactly one $n-1$ dimensional face which we call $F$, and 
furthermore $F$ intersects each of the coordinate axes. Note that each compact face of $N(H)$ is the intersection of $F$
with coordinate planes of various dimensions. [Sc] also shows that $H_F(x)$ is zero only at the origin. 
Setting various 
$x_i$ equal to zero, this means that for any compact face $F_i$ of $N(H)$, $H_{F_i}(x)$ must be nonvanishing
on $(\R - \{0\})^n$. In particular, by Lemma 4.1 the Hessian condition of Theorem 1.3 is satisfied. Hence by
Theorem 1.3, the conclusions of Theorem 1.1a) are valid if the Newton distance of $H$ is at most $2$, and 
the conclusions of Theorem 1.1b) are valid if it is greater than 2. Therefore $M$ is bounded 
on $L^p$ for $p > \max(2, d(H))$, and if $d(H) > 2$ this exponent is optimal whenever $y \notin T_y(S)$. 
This is equivalent to the theorem of [IoSa1] regarding the maximal operator $M$.

\noindent {\bf Example 3.} Consider the case when $n = 2$, as considered in [IkKeMu1] and [IkKeMu2]. In the 
latter paper, sharp $L^p$ boundedness for $M$ for $p > 2$ was proven in nearly full generality. Here, we get
sharpness in the setting of part b) of Theorem 1.1 as well as the corresponding part of Theorem 1.3. These 
situations correspond to $H(x_1,x_2)$ being in what are called 
adapted coordinates. Although there is always a coordinate change of the form $(x_1,x_2) \rightarrow 
(x_1, x_2 - a(x_1))$ or $(x_1,x_2) \rightarrow (x_1 - a(x_2),x_2)$ that places $H$ in adapted coordinates,
since $a$ may be nonlinear Theorem 1.1b) does not imply a real-analytic version of the general
result of [IkKeMu2]. There is one sense however, in which we get results beyond those of
[IkKeMu2]. Namely, the condition that $y \notin T_y(S)$ is required in the results of [IkKeMu2]. Since here we
only need this condition for sharpness, Theorem 1.1b) gives a new $L^p$ boundedness result for $M$ when $S$
is real-analytic, $y \in T_y(S)$, $H(x_1,x_2)$ is in
adapted coordinates, and $v(H)$ does not consist of exactly one vertex lying on a coordinate axis. 

It should
be pointed out that [IkKeMu2] also makes significant use of Newton polyhedrons, in two dimensions. Their
methods are substantially different from those used here however, in that they do not use damping techniques
or higher dimensional singularities methods.

\noindent {\bf Example 4.} Suppose $S$ has at least one nonvanishing principal curvature at $y$, as considered by 
Sogge in [So]. Then some second derivative of $H(x)$ is nonvanishing at the origin. After doing a linear 
change of variables, we can assume that in the Taylor expansion of $H(x)$ at the origin each $x_i^2$ appears. Thus
we may apply Theorem 1.3 to say the appropriate portion of Theorem 1.1 holds. Since the Newton distance is 
less than 2 and the zeroes of each $H_F(x)$ must have order at most 2, the conclusions of Theorem 1.1a) 
hold. As a result, whenever the Hessian condition is satisfied we recover the result of [So] that $M$ is 
bounded on all $L^p$ for $p > 2$. In the real-analytic case, the part of this theorem of [So] not covered is 
when the Hessian everywhere has rank 1, and in the polynomial case as indicated earlier this corresponds to 
when a linear coordinate change makes $H(x)$ a function of one variable. (Note however that these are 
nontrivial cases). 

Similarly, if $S$ is any smooth finite type surface, as long as the Taylor expansion of $H$ is nonvanishing at
the origin one can do a similar linear coordinate change to make $N(H)$ intersect each coordinate axis.
Once again we can apply Theorem 1.3 and say that one part of Theorem 1.1 holds (under the Hessian condition).
This gives the following extension of the theorem of Sogge-Stein [SoSt]. 

\noindent {\bf Corollary 1.4.} Suppose there are directions $u$ and $v$ such that the determinant of the
2 by 2 Hessian of $g$ in the
$u$ and $v$ directions does not vanish to infinite order at $z$. Then $M$ is bounded on $L^p$ for some
$p < \infty$.

Note that if $b$ denotes the order of the zero of $H(x)$ at the origin, then in Corollary 1.4 one may
take $p$ to be anything greater than $b$. We next turn to the surface measure Fourier transform analogues to Theorems 1.1 and 1.2. In the following 
theorems, $k$ denotes the dimension of the face (compact or not) of $N(H)$ that intersects the critical 
line $\{(t,...,t): t \in \R\}$ in its interior. If the line intersects $N(H)$ at a vertex we let $k = 0$.

\noindent {\bf Theorem 1.5.} Suppose $g(x)$ is either real-analytic or a smooth function such that the 
Newton polyhedron of $H(x)$ intersects each coordinate axis. Then there
is a neighborhood $U$ of $y$ such that if $\phi$ is supported in $U$ the following hold, where $T(\lambda)$
is as in $(1.3)$.

\noindent {\bf a)} If $d(H) < 2$, and each zero of each $H_F(x)$ on $(\R - \{0\})^n$ has order at most 2,
then there is a constant $C$ such that $|T(\lambda)| \leq C|\lambda|^{-{1 \over 2}}$ for $|\lambda| > 2$.

\noindent {\bf b)} If $d(H) \geq 2$ and each zero of each $H_F(x)$ on $(\R - \{0\})^n$ has order at most
$d(H)$, then there is a constant $C$ such that $|T(\lambda)| \leq C|\lambda|^{-{1 \over d(H)}}(\ln |\lambda|)^
{n - k}$ for $|\lambda| > 2$. If $d(H)$ is not an integer, the exponent $n - k$ can be improved to $n - k - 1$. 

\noindent {\bf c)} If the maximum order $m$ of any zero of any $H_F(x)$ on $(\R - \{0\})^n$ 
satisfies $m > \max(d(H),2)$ then there is a constant $C$ such that $|T(\lambda)| \leq C|\lambda|^
{-{1 \over m}}$ for $|\lambda| > 2$. 

The issue with the Hessian doesn't arise in Theorem 1.5 since this Hessian is used in order to be able to
apply a result of Sogge-Stein which is needed for the maximal operator only. Theorem 1.5 is sharp in the 
following ways.

\noindent {\bf Theorem 1.6.} Suppose $g(x)$ is real-analytic.

\noindent {\bf a)} If the hypotheses of Theorem 1.5 b) hold and $\phi$ is nonnegative with $\phi(y) \neq 0$ then 
$$\lim_{r \rightarrow \infty} \sup_{|\lambda| = r}{|T(\lambda)| \over |\lambda|^{-{1 \over d(H)}}(\ln |\lambda|)^
{n - k - 1}} > 0$$
\noindent {\bf b)} Suppose the hypotheses of Theorem 1.5 c) hold. Suppose there is a 
sequence  of points in $(\R - \{0\})^n$ with $a^k \rightarrow 0$ such that for each $k$, on a neighborhood of each
$a^k = (a^k_1,...,a^k_n)$ the function $H(x)$ is of the form 
$\xi_k(x)(x_i - b_k(x_1,...x_{i-1},x_{i+1},...,x_n))^m$ for some real-analytic $b_k$ and $\xi_k$ such that
$\xi_k(a^k) \neq 0$ and $b_k((a^k_1,...,a^k_{i-1},a^k_{i+1},...,a^k_n)) = a^k_i$. Then in any neighborhood $U$
of $y$ there is some $\phi$ supported in $U$ for which 
$$\lim_{r \rightarrow \infty} \sup_{|\lambda| = r}{|T(\lambda)| \over |\lambda|^{-{1 \over m}}} > 0$$

Because of the similarities between Theorems 1.5-1.6 and Theorems 1.1-1.3, the four examples given for
the maximal operator have direct analogues for the Fourier transforms of surface measures. The analogue
of example 4 is a somewhat trivial consequence of the Van der Corput lemma. The analogues of the other 
examples are not trivial. The analogue of example 1 says that when the Newton distance is at least 2,
the oscillatory integral estimates of [V] hold uniformly under linear perturbations of the phase. Example 2
gives sharp estimates for convex surfaces of finite-line type considered for example in [BrNaWa] and 
[CoDiMaMu]. Example 3 says that in adapted coordinates in two dimensions, the estimates of [V] for two-dimensional
oscillatory integrals hold uniformly under linear perturbations of the phase. Except for logarithmic factors, this also
follows from [IkKeMu2].

Lastly, we mention the following conjecture of Stein, who considered the $\alpha = {1 \over 2}$ case, and 
Iosevich-Sawyer, who extended it to all $\alpha \leq {1 \over 2}$:

\noindent {\bf Conjecture (Stein, Iosevich-Sawyer):} If $S$ is a smooth hypersurface and $\alpha \leq {1 \over 2}$ 
is such that $|T(\lambda)| < C_{\phi}(1 + |\lambda|)
^{-\alpha}$ for all $\phi$ supported on a sufficiently small neighborhood of $y$, then the maximal operator 
$M$ is bounded on $L^p$ for all $p > {1 \over \alpha}$. 

Observe that when $g(x)$ is real-analytic, Theorems 1.1a) and 1.5a) verify the conjecture under the
hypotheses of Theorem 1.1a).
Similarly, Theorems 1.5b)  and 1.6a) coupled with Theorem 1.1b) verify the conjecture when $d(H) > 2$ and each 
zero of each $H_F(x)$ on 
$(\R - \{0\})^n$ has order at most $d(H)$, assuming the vertex set of $N(H)$ does not consist solely of a
single vertex lying on a coordinate axis. The other results of this paper are consistent with the conjecture.
However, we do not make any guesses here about whether or not the full
conjecture holds, as in this paper we work only with a coarse resolution of singularities and there 
may be additional issues arising when the finer aspects of the singularities of $H(x)$ become pertinent. 

\noindent {\bf 2. Outline of method.}

In several papers such as  [IoSa1] [IoSa2] [NaSeWa] [SoSt] one successful technique that has helped in 
proving $L^p$ bounds for $M$ has been to interpolate
between $L^2$ and $L^{\infty}$ boundedness of damped versions of $M$. Specifically, for a properly selected
function $p(x)$ on $\R^{n+1}$ one defines the measure $\sigma_z$ by 
$d\sigma_z = |p(x)|^z d\sigma$. Since we will be working in the situation where $S$ is the graph
of a function $g(x)$ in the first $n$ variables, for our purposes $p(x)$ will be a function of the 
first $n$ variables. Often $p(x)$ is defined directly in terms of $g(x)$ itself. The idea is that if $p(x)$ 
is chosen properly and $M_z$ denotes the maximal operator with $\sigma$ 
replaced by $\sigma_z$, then for some $b > 0$, $M_z$ is bounded on $L^2$ for $Re(z) > b$
and for some $a > 0$, $M_z$ is bounded on $L^{\infty}$ for $Re(z) > -a$. For a small $\epsilon > 0$ one then
uses complex interpolation (the variant for maximal operators, that is) on $M_z$ between the lines 
$Re(z) = -a + \epsilon$ and $Re(z) = b + \epsilon$ to conclude that $M = M_0$ is bounded on 
$L^{p(\epsilon)}$ for a certain value $p(\epsilon)$. Letting $\epsilon$ go to zero, $p(\epsilon)$ will 
converge to some $p_0$ such that $M$ is bounded on $L^p$ for all 
$p > p_0$. If all goes well, $p_0$ is the optimal value. There are some variations on this method; for 
example sometimes one needs to define $\sigma_z = 
e^{z^2}|p(x)|^z \sigma$ to retain boundedness of $M_z$ as $|Im(z)| \rightarrow \infty$. 

Ever since the earliest papers in this area such as [St1], one way of proving $L^2$ boundedness of maximal
operators of the kind considered in this paper has been to reduce the
problem to proving decay estimates for the Fourier transform of an associated surface measure. The version 
of this idea we will use here is the following consequence of [SoSt].

\noindent {\bf Theorem ([SoSt])} Suppose for some $z$ there are $C, \epsilon > 0$ such that for all
multiindices $\alpha$ with $|\alpha| = 0,1$ the measure $d\mu_z =  |p(x)|^z d\sigma$ satisfies  
$$|\partial^{\alpha} \hat{\mu}_z (\lambda)| < C(1 + |\lambda|)^{-{1 \over 2} - \epsilon} \eqno (2.1)$$
Then there is a constant $C'$ depending on $C$ and $\epsilon$ such that $||M_zf||_2 \leq C'||f||_2$ 
for all $f \in L^2$. 

In practice, if one has $(2.1)$ for $\alpha = 0$, it will generally automatically hold for all the first 
derivatives since the effect of taking such a derivative is to replace the cutoff function by another one. In
fact, one sometimes gets better behavior for the first derivatives due to the new cutoff functions being 
zero at the origin. 
To motivate the choice of $p(x)$ in this paper, we examine $\hat{\mu}_z(\lambda)$. We will always choose $p(x)$ 
to be a function of the first $n$ variables only. Then the function $\hat{\mu}_z(\lambda)$ is given by 
$$\hat{\mu}_z(\lambda) = \int e^{-i\lambda_1 x_1 - ... - i\lambda_n x_n - i\lambda_{n+1}g(x_1,...,x_n)}
|p(x_1,...,x_n)|^z \psi(x_1,...,x_n)\,dx_1...\,dx_n \eqno (2.2)$$
Since invertible linear transformations in $\lambda$ don't affect the truth of $(2.1)$, as in the last
section we can replace $g(x)$ by $H(x) = g(x + z) - g(z) - \nabla g(z) \cdot x$ without changing whether or
not $(2.1)$ holds. So letting $P(x_1,...,x_n) = p(x_1 + z_1,...,x_n + z_n)$ and $\Psi(x_1,...,x_n) = 
\psi(x_1 + z_1,...,x_n + z_n)$ we must show $G_z(\lambda)$ and its first partials are bounded by the 
right-hand side of $(2.1)$, where 
$$G_z(\lambda) = \int e^{-i\lambda_1 x_1 - ... - i\lambda_n x_n - i\lambda_{n+1}H(x_1,...,x_n)}
|P(x_1,...,x_n)|^z \Psi(x_1,...,x_n)\,dx_1...\,dx_n \eqno (2.3)$$
As might be expected from the fact that the statement
of Theorem 1.1 involves the Newton polyhedron of $H$ in a way analogous to the theorems of [G1], the choice 
of $P(x_1,...,x_n)$ and various other aspects of the proof of our results will draw from [G1].
In a sense, we will consider 
the phase in $(2.3)$ as a linear perturbation of $\lambda_{n+1}H(x_1,...,x_n)$, the one-parameter phase 
traditionally considered for oscillatory integrals, and we will extend the methods of [G1] to cover such 
linear perturbations and the damping/interpolation methods needed for the maximal operators. It should
be pointed out that the relevance of singularities and the Newton polyhedron in studying oscillatory 
integrals has been known
for some time, such as in [V] for scalar oscillatory integrals and in [PhSt] for operator versions.

To see how our arguments will proceed, we first consider the situations where each $H_F(x)$ has no zero of
order greater than 2 in $(\R - \{0\})^n$. This covers part a) and some of part b) of 
Theorem 1.1. Note that these situations include examples 1 and 2 of the last section. As in the analysis 
of [G1], our first step is to divide the domain of $(2.3)$ into $2^n$ octants via the $n$ coordinate 
hyperplanes $x_j = 0$. We focus our attention on the octant where each $x_j > 0$ as the other ones are dealt
with in the analogous fashion. Similar to in [G1] (and also [PhSt] incidentally), we divide this 
octant into rectangles of the form $R = \prod_{i=1}^n [2^{-k_i-1},2^{-{k_i}}]$.  Letting the damping
factor $P(x)$ just be 1 for now, the portion of $(2.3)$ corresponding to $R$ is given by
$$I_R = \int_R e^{-i\lambda_1 x_1 - ... - i\lambda_n x_n - i\lambda_{n+1}H(x_1,...,x_n)} \Psi(x_1,...,x_n)\,dx_1
...\,dx_n \eqno (2.4)$$
Let $v(H)$ denote the set of all vertices of $N(H)$, and define $H^*(x_1,...,x_n) = \sum_{v \in v(H)} |x|^
v$. We will see that the fact that each $H_F(x)$ has no zero of order greater than 2 
in $(\R - \{0\})^n$ enables one to divide $R$ into boundedly many rectangles $R_j$ such that on each $R_j$ 
there is a derivative  $\partial_v = \sum_i a_i \partial_{x_i}$ with each $|a_i| < 2^{-{k_i}}$ such that 
$$|\partial_v^2 H(x_1,...,x_n)| > C_0 H^*(x_1,...,x_n) \eqno (2.5)$$
Here $C_0$ is a constant depending on the function $H$. This will allow us to use a one-dimensional Van der 
Corput lemma on each $R_j$. Once the resulting bounds are added over all $j$, we will obtain the estimate 
$$ |I_R| \leq C_1 |R| (|H^*(R)| |\lambda|)^{-{1 \over 2}} \eqno (2.6)$$
Here $C_1$ again denotes a constant depending only on $H$, $|R|$ is the measure of the rectangle $R$, and 
$|H^*(R)|$ is shorthand for the supremum of $H^*(x_1,...,x_n)$ over $R$. Note that $H^*(x_1,...,x_n)$ only 
varies by a factor of $C$ on any given rectangle, so up to a constant $H^*(x_1,...,x_n)$ is equal to 
$|H^*(R)|$ everywhere on $R$. 
Adding $(2.6)$ over all $R$ will then give an upper bound for the overall integral $(2.3)$ for $G_z(\lambda)$
when $P(x)$ is 1. 

Suppose now one chooses $P(x)$ to be a nonconstant function, but still a function that
varies by at most a constant on each $R$. Then we get the following analogue of $(2.6)$:
$$|I_R| \leq C_2 |R||P^*(R)|^{Re(z)}  (|H^*(R)| |\lambda|)^{-{1 \over 2}} \eqno (2.7)$$
Actually, in our arguments $C_2$ will increase linearly in $|Im(z)|$, but the above-mentioned 
replacement $M_z$ by $e^{z^2}M_z$ in our interpolations will always take care of such factors. 
The question now is how to choose $P(x)$. In some papers such as [IoSa1] - [IoSa2] one lets the damping 
factor be $H(x)$ itself. In this paper the zeroes of $H(x)$ make $H(x)$ unsuitable for this purpose, but by 
choosing $P(x) = H^*(x)$ one can eliminate the effect of all zeroes except for the main one at the 
origin which determines the overall integrability behavior of negative powers of $|H|$.
So we preliminarily select $P(x) = H^*(x)$. 

Suppose for the time being that $d(H) \geq 2$, in addition to the above assumption that any zero of any
$H_F(x)$ has order at most 2. By [G1], the function $H^*(x)^{-t}$ is integrable for $t < 
{1 \over d(H)} $, so by $(2.7)$ $|\sum_R I_R|$ is finite for $Re(z) > {1 \over 2} - {1 \over d(H)}$; in fact
for such $z$ we have $|\sum_R I_R| \leq C_z|\lambda|^{-{1 \over 2}}$. In order to apply the theorem of Sogge and
Stein however, we need an extra epsilon of decay in the exponent. For this we tack on an additional damping 
factor, and this is where the assumed condition on the 
2 by 2 Hessian determinants comes in. Specifically, we select directions $u$ and $v$ such that the determinant
$D(x)$ of the Hessian of $H$ in the $u$ and $v$ directions does not vanish to infinite order at the origin.
Instead of choosing $P(x) = H^*(x)$, we choose $P(x) = |D(x)|^{\delta}H^*(x)$. This will give the added 
$\epsilon$ needed for $(2.1)$ hold in the lemma of Sogge-Stein.

The way the adjusted damping function will accomplish this is as follows. We will write the function 
$G_z(\lambda)$ of $(2.3)$ as $G_z^1(\lambda) + G_z^2(\lambda)$, where $G_z^1(\lambda)$ is the integral over the set 
where
$|D(x)| < |\lambda|^{-\eta}$, and $G_z^2(\lambda)$ is the integral over the set where $|D(x)| \geq 
|\lambda|^{-\eta}$. Here $\eta$ is a small constant depending only on the dimension $n$. If one argues as above
in the first integral, the additional factor of $|D(x)|^{\delta}$  gives an additional factor of 
$C|\lambda|^{-\eta \delta}$ and the above analysis yields the estimate $|G_z^1(\lambda)| < 
C|\lambda|^{-{1 \over 2} - {\eta \delta} Re(z)}$. This improves the exponent to satisfy $(2.1)$.
For $G_z^2(\lambda)$, the fact that the Hessian determinant in the $u$ and $v$ directions is bounded below
by the relatively large factor $C'|\lambda|^{-\eta}$ will allow us to argue similarly to the 2-dimensional
nondegenerate case in the $u$-$v$ directions, and then integrate the result. Because the 2-d nondegenerate
case gives an estimate of $C|\lambda|^{-1}$, there will be a lot of slack in this portion of the argument,
and we will get the relatively strong decay of $C|\lambda|^{-{3 \over 5}}$ for $|G_z^2(\lambda)|$. Adding the
upper bounds for $|G_z^1(\lambda)|$ and $|G_z^2(\lambda)|$ will give us the needed bounds for $|G_z(\lambda)|$.

As for proving the $L^{\infty}$ boundedness of $M_z$ for $z$ satisfying $Re(z) > -a$ for some $a > 0$, as is 
typical in such arguments $a$ will be maximal such that the measure $|\mu_z|$ is finite for $Re(z) > -a$. 
Thus $a$ is the supremum of all $t > 0$ such that $(|D(x)|^{\delta}H^*(x))^{-t}$ is integrable. One
interpolates this $L^{\infty}$ boundedness of $M_z$ with the above $L^2$ boundedness to obtain the 
relevant part of Theorem 1.1
as follows. Fix $\epsilon > 0$. The supremum $a'$ of the $t$ which $H^*(x)^{-t}$ is integrable can be 
obtained using the results of [G1]. If one lets $\delta$ go to zero, $a$ will converge to $a'$. So if 
$\delta$ is small enough, for any $z$ with $Re(z) = -a' + \epsilon$, $M_z$ is bounded on  
$L^{\infty}$ with constants uniform in $|Im(z)|$. If one interpolates this with the above $L^2$ boundedness of $M_z$ for $Re(z) = {1 \over 2} 
- {1 \over d(H)} + \epsilon$, we obtain a value of $p'$ for which $M$ is bounded on $L^p$ for $p > p'$. One then 
lets $\epsilon$ go to zero, and $p'$ will converge to some $p_0$ for which we will have
shown that $M$ is bounded on $L^p$ for $p > p_0$. This $p_0$ will be the exponent given by Theorem 1.1.

Suppose now that some $H_F(x)$ has a zero of order greater than 2 on $(\R - \{0\})^n$.
Then the above arguments break down; in particular $(2.5)$ 
no longer holds. However, one still has that the left-hand side of $(2.5)$ is bounded above by $CH^*(x)$. In
fact we will show that
$$\sum_{i,j} 2^{-k_i -k_j} |\partial_{x_ix_j}^2 H(x)| \leq CH^*(x)$$
To deal with the new situation we introduce an additional damping factor. To understand what 
is needed, we first divide the integral of $(2.4)$ into pieces on each of which $\sum_{i,j} 2^{-k_i -k_j}
|\partial_{x_ix_j}^2 H(x)|$ is between $2^{-l}H^*(x)$ and $2^{-l+1}H^*(x)$ for some $l$. We correspondingly 
write $I_R = \sum_l I_R^l$. The effect of shrinking $\sum_{i,j} 2^{-k_i -k_j} |\partial_{x_ix_j}^2 H(x)|$
by a factor of $2^l$ from what it would be if $(2.5)$ held is to cause $I_R^l$ to be $C2^{{l \over 2}}$
times the estimate given by $(2.6)$. Thus to counteract this effect, we would like to introduce an
additional damping factor of $\big(\sum_{i,j} 2^{-k_i -k_j}|\partial_{x_ix_j}^2 H(x)|\big)^{1 \over 2}$
into the $L^2$ estimates. For technical reasons, it is more convenient to choose as our damping factor 
$H^{**}(x) = \big(\sum_{i,j} x_i^2 x_j^2 (\partial_{x_ix_j}^2 H(x))^2\big)^{1 \over 4}$. 

The presence of $H^{**}(x)^z$ changes the factor of $H^*(x)^z$ needed in the damping factor. One way of
seeing how this works is the following. If we were back into the situation where all zeroes of all 
$H_F(x)$ are of order 2 or less, $H^{**}(x)$ would be comparable in magnitude to $H^*(x)^{1 \over 2}$.
Even when the zeroes are not all of order 2 or less, over most of a given rectangle $H^{**}(x)$ is
still comparable in magnitude to $H^*(x)^{1 \over 2}$. So in general, if 
instead of $H^{**}(x)$ we insert the damping factor $H^*(x)^{-{1 \over 2}}H^{**}(x)$ into the $L^2$ 
estimates, this new factor serves to eliminate the effect of zeroes of the $H(x)$ without any additional 
impacts. We want this factor to be there at the "edge" of when the $L^2$ estimates hold, namely for 
$Re(z) = {1 \over 2} - {1 \over d}$, where $d$ is shorthand for $d(H)$. Hence the original damping factor of
$H^*(x)^z$ gets replaced by the damping factor $H^*(x)^z\big(H^*(x)^{-{1 \over 2}}H^{**}(x)
\big)^{z \over {1 \over 2} - {1 \over d}}$ $ = H^*(x)^{-{2 \over d - 2}z} H^{**}(x)^
{{2d \over d - 2}z}$. It will actually be more convenient to use the ${d - 2 \over 2d}$ power of this.  
Since we still need the $|D(x)|^{\delta z}$ factor for the same reasons as before, we obtain
an overall damping factor of 
$$P(x)^z = \big(|D(x)|^{\delta}H^*(x)^{-{1 \over d}} H^{**}(x)\big)^z \eqno
(2.8)$$
Our analysis will then allow us to prove the estimate $(2.1)$ needed for the Sogge-Stein lemma when $Re(z) 
> 1$. This will give the $L^2$ estimates needed for such $M_z$. For $L^{\infty}$
boundedness of $M_z$, like above $M_z$ is bounded on $L^{\infty}$ for $Re(z) > -a$, where $a$ is the 
supremum of the $t$ for which 
$P(x)^{-t}$ is integrable. Once again, as $\delta$ tends to zero $a$ converges to the 
analogous quantity $a'$ for $H^*(x)^{-{1 \over d}} H^{**}(x)$. The 
results of [G1] will allow us to determine $a'$. Then for any $z$
with $Re(z) > a'$, if $\delta$ is sufficiently small one has $L^{\infty}$ boundedness for $M_z$. 
Interpolating the $L^2$ and $L^{\infty}$ boundedness for $M_z$ will give the rest of Theorem 1.1.

Although the above assumed that $d > 2$ and some $H_F(x)$ has a zero in $(\R - \{0\})^n$ of order greater
than 2, the above damping factor will still work when $d > 2$ and any zero of any $H_F(x)$ has order 2 or
less. The reason is that in this case $H^{**}(x)$ is everywhere comparable in magnitude to $H^*(x)^{1 
\over 2}$, so the damping factor serves the same purpose as a damping factor $|D(x)|^{\delta}$ times
a power of $H^*(x)$, and we saw above that such a factor is appropriate for this situation.

When $d \leq 2$ and some $H_F(x)$ has a zero in $(\R - \{0\})^n$ of order greater than 2, one no longer
needs to include the original damping factor of $H^*(x)$, just the second damping factor $H^*(x)^
{-{1 \over 2}}H^{**}(x)$. The $|D(x)|^{\delta}$ factor is still needed for the same reason as before.
Hence one can merge the $d \leq  2$ and $d > 2$ cases by writing the damping factor in all cases as
$$P(x)^z = \big(|D(x)|^{\delta}H^*(x)^{-{1 \over \max(d,2)}} H^{**}(x)\big)^z \eqno
(2.9)$$
We still haven't considered the case where $d \leq 2$ and any zero of any $H_F(x)$ has order 2 or less 
(the situation of Theorem 1.1a). But here we may still use damping factor $(2.9)$. For in this case, the fact 
that $H^{**}(x) \sim H^*(x)^{-{1 \over 2}}$ means that
$P(x) \sim |D(x)|^{\delta}$. If $\delta$ is small enough, then $M_z$ is bounded on $L^{\infty}$ on 
$Re(z) = -1$ with uniform constants. On the other hand we will see that for any $b > 0$, if $Re(z)$ then
by adding up the
damped version of $(2.7)$ over all $R$ one has that $|G_z(\lambda)| < C_z|\lambda|^{-{1 \over 2} - \epsilon}$. 
So by the Sogge-Stein lemma, $M_z$ is bounded on $L^2$ on whenever $Re(z) = b$. Once again $C_z$ grows
linearly in $|Im(z)|$ and thus is controllable by the additional $e^{z^2}$ factor
in the interpolation. Interpolating this with the above $L^{\infty}$ bounds for a sequence of $b$'s 
tending to zero gives $L^p$ boundedness of $M$ whenever $p > 2$, implying part a) of Theorem 1.1. 

Proving the estimates for Fourier transforms of surface-supported measures given by Theorem 1.5 will be
substantially simpler than the above since we won't have to worry about damping factors and interpolation.
Once again we will divide the domain of the integral into the rectangles $R$. This time, each $R$ will be
subdivided into boundedly many subrectangles $R_j$ on which we will be able to use the methods of [G1]
and directly apply the appropriate Van der Corput lemma to get the necessary estimate for 
the integral over $R_j$. Adding the resulting estimates over all rectangles will give Theorem 1.5. 
Actually, for the purposes of these arguments the terms  $-i\lambda_1 x_1... - i\lambda_n x_n$ appearing in the phase of $(1.7)$ can be 
viewed as a linear perturbation of the phase $i\lambda_{n+1} H(x_1,...,x_n)$ which do not affect the 
applicability of the arguments of [G1] in any major way; in our proofs one always uses a Van der Corput 
lemma for functions with a nonvanishing derivative of order two or higher. 

\noindent {\bf 3. Preliminary lemmas and proofs of Theorems 1.5-1.6.}

Suppose $f(x)$ is a real-analytic function on a neighborhood of the origin in $\R^n$ such that $f(0) 
= 0$. We now describe
the decomposition of a small neighborhood $U$ of the origin done in [G1] that we will use in this paper. This
decomposition is done according to the Newton polyhedron $N(f)$ of $f(x)$ as follows. First, one divides
the neighborhood into $2^n$ octants via the coordinate hyperplanes $x_k = 0$. We describe the decomposition
for the open octant where each $x_k > 0$, as the other ones are given by reflecting the decomposition about 
varous hyperplanes $x_k = 0$. 

Let $F_{i1},...,F_{ik_i}$ denote the set of compact faces of $N(f)$ of dimension $i$. 
As in previous sections, we consider a vertex of $N(f)$ to be a face of dimension zero. For each $F_{ij}$
we will have a finite collection of open wedges $W_{ijp}$ whose closures each contains the origin.
Each $W_{ijp}$ is defined through some monomial inequalities $\{x \in U: x^{\alpha_1} < Cx^{\alpha_2}\}$ or 
$\{x \in U: x^{\alpha_1} < Cx^{\alpha_2}\}$. Up to a set of measure zero, the union of all the $W_{ijp}$
is the whole neighborhood $U$ and furthermore, the decomposition can be done so that the $W_{ijp}$ has a
number of properties which we will give as Theorems 3.1-3.3. The first is a direct consequence of Theorem 
2.0 of [G1] and also of Theorem 3.2 of [G2]:

\noindent {\bf Theorem 3.1.} Let $v(f)$ denote the set of vertices of $N(f)$. There are $A_1, A_2 > 1$ 
depending on the function $f(x)$ such that if $C_0,...,C_n$
are constants with $C_0 > A_1$ and $C_{i+1} > C_i^{A_2}$ for all $i$, then one can define the $W_{ijp}$ 
so that

\noindent {\bf a)} Let $i < n$. If the following two statements hold, then  $x \in \cup_p W_{ijp}$.

\noindent {\bf 1)} If $v \in v(f) \cap F_{ij}$ and  $v' \in v(f) \cap
(F_{ij})^c$ we have $ x^{v'} < C_n^{-1} x^v$. \parskip = 0pt

\noindent {\bf 2)} For all $v, w \in v(f) \cap F_{ij}$ we have $C_i^{-1}x^w < x^v < C_ix^w$. \parskip = 12pt

\noindent {\bf b)} There is a $\delta > 0$ depending on $N(f)$, and not on $A_1$ or $A_2$, such that if 
$x \in \cup_p W_{ijp}$, then the following two statements hold.

\noindent {\bf 1)} If $v \in v(f) \cap F_{ij}$ and  $v' \in v(f) \cap
(F_{ij})^c$ we have $ x^{v'} < C_{i+1}^{-\delta} x^v$. \parskip = 0pt

\noindent {\bf 2)} For all $v, w \in v(f) \cap F_{ij}$ we have $C_i^{-1}x^w < x^v < C_ix^w$. \parskip = 12pt

Informally, Theorem 3.1 gives a way of saying that the vertices of $F_{ij}$ dominate the Taylor series
of $f$ when $x$ is in one of the $W_{ijp}$. Another way of making this precise which we will make significant
use of is the following lemma, called Lemma 2.1 in [G1]:

\noindent {\bf Theorem 3.2.} Suppose $x \in \cup_p W_{ijp}$. Let $V \in v(f)$ be such that $x^V \geq x^v$ 
for all $v \in v(f)$; if there is more than one such vertex let $V$ be any of them. Then if
$A_1$ is sufficiently large and $\eta$ is sufficiently small, for any positive $d$ one has the following estimate:
$$\sum_{\alpha \notin F_{ij}} |s_{\alpha}| |\alpha|^d x^{\alpha} < K(C_{i+1})^{-\delta''}x^V$$
Here $K$ is a constant depending on $d$ as well as the function $f(x)$, and $\delta'' > 0$ is a constant
depending on the Newton polyhedron of $f$. 

Note that in Theorem 3.1, one can increase $A_1$ without affecting the truth of the theorem, so stipulating
that $A_1$ must be sufficiently large in Theorem 3.2 is entirely consistent with Theorem 3.1. Next we give
the following straightforward corollary of Theorem 3.2 (Corollary 2.2 of [G2]):

\noindent {\bf Corollary.} There is a constant $C$ depending on the function $f(x)$ such that on a 
sufficiently small neighborhood of the origin $|f(x)| \leq C \sum_{v \in v(f)} x^v$.

The next theorem (Theorem 2.2 of [G1]) describes the monomial map that one can do on a given $W_{ijp}$ that 
converts it to a set
comparable to a cube on which each $x^v$ for $v \in F_{ij}$ becomes the same monomial $x^{v'}$, and on which
each $x^v$ for $v \notin F_{ij}$ becomes a monomial $x^w$ (depending on $v$) for which each $w_i \geq (v')_i$ 
with at least one inequality strict. Once again, it assumes that $A_1$ and $A_2$ are sufficiently large,
and as before this is consistent with the previous theorems.

\noindent {\bf Theorem 3.3.} If $A_1$ and $A_2$ are sufficiently large, to each $W_{ijp}$ 
there is a bijective map $\beta_{ijp}: Z_{ijp}  \rightarrow W_{ijp}$ depending
on $N(f)$ and $(i,j,p)$ such that each component of 
$\beta_{ijp}(z)$ is a monomial in $(z_1^{{1 \over N}},...,z_n^{{1 \over N}})$ for some positive integer $N$,
and such that for some $\mu ' > 0$ we have 

\noindent {\bf a)} When  $i = 0$, $(0, \mu ')^n \subset Z_{ijp} \subset (0,1)^n$.\parskip = 3 pt

\noindent {\bf b)} When $i > 0$, there are sets $D_{ij} \subset (C_i^{-e}, C_i^e)^i$ for some $e > 0$ 
depending on $N(f)$ such that $(0, \mu ')^{n-i} \times D_{ij} \subset Z_{ijp} \subset 
(0, 1)^{n-i} \times D_{ij}$

\noindent {\bf c)} When $i > 0$, write $z \in \R^n$ as $(\sigma,t)$ where
$\sigma \in \R^{n-i}$ and $t \in \R^i$. For any $v \in N(f)$, denote by $\sigma^{v'}t^{v''}$ the
function in $z$ coordinates that $x^v$ transforms into under the $x$ to $z$ coordinate change. When
$i = 0$, write $z = \sigma$ and for $v \in N(f)$ denote by $\sigma^{v'}$ the the function
$x^v$ transforms into. 
Then for any $v_1,v_2 \in F_{ij}$ we have $v_1' = v_2'$, while if $v_1 \in F_{ij}$
and $v_2$ is in $N(f)$ but not in $F_{ij}$, then $(v_2')_k \geq (v_1')_k$ for all $k$ with at least 
one component strictly greater. \parskip = 12 pt

\noindent We need one more lemma from [G1]:

\noindent {\bf Lemma 3.4 (Lemma 3.1 of [G1]).} Suppose $m_1,...,m_n$ are nonnegative numbers not all zero. 
Let $M = \max_i
m_i$, and let $l$ denote the number of $m_i$ equal to $M$. Then if $|E|$ denotes Lebesgue measure,
we have the following for all $0 < \delta < 1$, where $C$ and $C'$ are constants depending on the 
$m_i$. 

\noindent {\bf a}) \centerline{$C |\ln \delta|^{l-1}\delta^{1 \over M} < |\{x \in (0,1)^n: x_1^{m_1}...
x_n^{m_n} < \delta\}| < C' |\ln \delta|^{l-1}\delta^{1 \over M}$}

\noindent {\bf b)} If $M < 1$, then 
$$C \delta  < \int_{\{x \in (0,1)^n: {\delta \over x_1^{m_1}...x_n^{m_n}} < 1\}} {\delta \over
x_1^{m_1}...x_n^{m_n}}\,dx < C' \delta $$
\noindent {\bf c)} If $M = 1$, then 
$$C |\ln \delta|^l \delta  < \int_{\{x \in (0,1)^n:  {\delta \over x_1^{m_1}...x_n^{m_n}} < 1\}} {\delta \over x_1^{m_1}...x_n^{m_n}}\,dx
< C' |\ln \delta|^l \delta $$
\noindent {\bf d)} If $M > 1$, then 
$$\int_{\{x \in (0,1)^n: {\delta \over x_1^{m_1}...x_n^{m_n}} < 1\}} {\delta \over x_1^{m_1}...x_n^{m_n}}\,dx < 
C' |\{x \in (0,1)^n: x_1^{m_1}...x_n^{m_n} < \delta\}|$$ 
\noindent The next lemma says that if one of the polynomials $f_F$ has a zero of order $\leq 2$ at some $y \in 
(\R - \{0\})^n$, then in fact some second partial $\partial_{x_i x_j}^2 f_F$ is nonzero at $y$. 

\noindent {\bf Lemma 3.5.} Suppose $f(x)$ is a smooth function satisfying $f(0) = 0$ and $\nabla f(0) = 0$
which has a nonvanishing
Taylor expansion at the origin. Let $f_F(x)$ be one of the polynomials of Definition 1.2. Suppose $y \in
(\R - \{0\})^n$ is such that $f_F(y)$ is nonzero or has a zero of order 1 at $y$. Then there are variables 
$x_i$ and $x_j$ such that $\partial_{x_i x_j}^2 f_F(y) \neq 0$.

\noindent {\bf Proof.} First consider the case where $f_F(y) \neq 0$. Since $F$ is a face of $N(f)$, we
may let $\alpha_1,...\alpha_n, a$ be positive numbers such that $\alpha \cdot v = a$ for all vertices
$v$ of $N(f)$ on $F$, and $\alpha \cdot v > a$ for all vertices not on $F$. Then as a function of $t$, 
$f_F(y_1t^{\alpha_1},...,y_nt^{\alpha_n})$ is of the form $C(y)t^a$ with $C(y) \neq 0$. Because $f(0) = 
0$, also $f_F(0) = 0$ and therefore $a > 0$. Thus if one 
takes the derivative with respect to $t$
of the equation $f_F(y_1t^{\alpha_1},...,y_nt^{\alpha_n}) = C(y)t^a$ and sets $t = 1$,
one gets $\sum_i \alpha_i y_i \partial_{x_i} f_F(y)$ is equal to the nonzero quantity $C(y)a$. Hence at least
one term $y_i \partial_{x_i} f_F(y)$ is nonzero. Since $y_i \neq 0$, this means some $\partial_{x_i} f_F(y)$
is nonzero and we have reduced to the case where some first partial of $f_F$ is nonzero at $y$. 

So to prove the lemma, we must show that if some first partial is nonzero at $y$, then so is some second
partial. To do this, we perform the above argument on $\partial_{x_i} f_F(y)$. This time, for exponents 
$v$ in the nonvanishing terms of $\partial_{x_i} f$'s Taylor expansion, the quantity $\alpha \cdot v$ is 
minimized for the monomials appearing in $\partial_{x_i} f_F$. Furthermore, since $f(x)$ has a critical
point at the origin, $\partial_{x_i} f_F$ still has a zero at the origin. Hence the above argument applies
once again, and we get a nonvanishing second partial as needed. 

\noindent {\bf Lemma 3.6.} Suppose $f(x)$ is a smooth function satisfying $f(0) = 0$ and $\nabla f(0) = 0$
which has a nonvanishing
Taylor expansion at the origin. Suppose that either $f(x)$ is real-analytic or a smooth function whose Newton polyhedron intersects
each coordinate axis. Let $v(f)$ denote the set of vertices of $N(f)$ and define $f^*(x) = \sum_{v \in v(f)}
|x|^v$. Let $m$ denote the maximum order of any zero of any $f_F(x)$, and let $M =
\max(2,m)$. Then there is a $\delta > 0$ and a neighborhood $U$ of the origin such that for each $x \in U$
there is a multiindex $\alpha$ with $0 \leq |\alpha| \leq m$ with
$$|x^{\alpha}\partial^{\alpha} f(x)| \geq \delta f^*(x) \eqno (3.1)$$
One may also arrange it so that $2 \leq |\alpha| \leq M$ for all $x$.

\noindent {\bf Proof.} We first suppose $f(x)$ is real-analytic. We apply the resolution of singularities 
algorithm giving Theorems 3.1-3.3 to $f(x)$, obtaining the resulting sets $W_{ijp}$. 
Suppose $F$ is any face of $N(f)$ of some dimension $i \geq 0$. Let $j$ be such that the sets
$W_{ij1},...,W_{ijk_{ij}}$ correspond to the face $F$, and let $Z_{ijp}$, $\sigma$ and $t$ be as in Theorem 
3.3. 

Let $\alpha$ be any multiindex such that $x^{\alpha}\partial^{\alpha} f_F$ is not the 
zero polynomial. Note that the exponent $v$ of any 
term $cx^v$ of the polynomial $x^{\alpha}\partial^{\alpha} f_F(x)$ is on $F$. Therefore if one writes 
$cx^v$ as $c\sigma^{v'}t^{v''}$ in the $z$ coordinates, by part c) of Theorem 3.3 each component
$(v')_k$ of $v'$ has the minimal possible value for any term in $f$'s Taylor expansion. So in the
$z$ coordinates, $x^{\alpha}\partial^{\alpha} f_F(x)$ becomes $\sigma^{v'}r_{\alpha,F}(t)$, where 
$r_{\alpha,F}(t)$
is some rational function of $t$. If $x \in (\R - \{0\})^n$ with $x^{\alpha}\partial^{\alpha}f_F(x)
\neq 0$, then if $(\sigma,t)$ are the $z$-coordinates of $x$ we have $r_{\alpha,F}(t) \neq 0$. (If $i = 0$,
then $x^{\alpha}\partial^{\alpha} f_F(x)$ is simply of the form $c\sigma^{v'}$ in the $z$ coordinates,
where $c \neq 0$).

By our assumptions, for any $x \in (\R - \{0\})^n$ there is a multiindex $\alpha$ of order $\leq m$ such 
that $x^{\alpha}\partial^{\alpha} f_F(x) \neq 0$. By Lemma 3.5, one can also take $\alpha$ to satisfy $2 
\leq \alpha \leq M$. Hence for a given $F$ of dimension $i > 0$, for all $t \in (\R - \{0\})^i$, there's some such $\alpha$ for 
which $r_{\alpha,F}(t) \neq 0$. Since by Theorem 3.3 the $t$ in $Z_{ijp}$ are restricted to a compact set
not intersecting the coordinate hyperplanes, there is some $\epsilon_0 > 0$ such that for all $(\sigma,t)$ 
in $Z_{ijp}$, for some such $\alpha$ we have
$$|r_{\alpha,F}(t)| > \epsilon_0 \eqno (3.2)$$ 
This in turn implies that $|\sigma^{v'}r_{\alpha,F}(t)| > \epsilon_0 \sigma^{v'}$. Translating back into the
original coordinates now, we conclude that if $v$ denotes some vertex of $N(f)$ on $F$, then if $x \in 
 W_{ijp}$ there is some $\alpha$ with $0 \leq |\alpha| \leq m$ such that
$$|x^{\alpha}\partial^{\alpha} f_F(x)| > \epsilon_0 x^v \eqno (3.3)$$
One may also take $\alpha$ to satisfy $2 \leq |\alpha| \leq M$. Although we assumed $i > 0$ in the derivation
of $(3.3)$, observe that it is trivially true when $i = 0$ as well. For any $i$, let $C_i$ be as in Theorem 3.1.
By Theorem 3.1b), if $x \in W_{ijp}$ then $x^v \geq C_i^{-1}x^w$ for all other vertices $w$ of $N(f)$. 
Thus we can replace $(3.3)$ by
$$|x^{\alpha}\partial^{\alpha} f_F(x)| > \epsilon_1 f^*(x) \eqno (3.3')$$
It is important to note that $\epsilon_1$ depends on $C_i$ but not $C_j$ for $j > i$; we will now
need to increase such $C_j$ while keeping $\epsilon_1$ fixed. Namely, let $\sum_{\gamma}f_{\gamma}x^{\gamma}$
denote the Taylor series of $f(x)$ about the origin. If $f(x)$ is real-analytic, then we have
$$|x^{\alpha}\partial^{\alpha}f(x) - x^{\alpha}\partial^{\alpha}f_F(x)| \leq \sum_{\gamma \notin F}
|\gamma|^m |f_{\gamma}|x^{\gamma} \eqno (3.4)$$ 
(If instead of $0 \leq |\alpha| \leq m$ we have $2 \leq |\alpha| \leq M$, then $m$ is replaced 
by $M$ in $(3.4)$). 
By Theorem 3.2, by making $C_{i+1}$ is sufficiently large for a fixed $C_i$ (i.e. one first chooses $C_0$,
then $C_1$, and so on), one can make the difference
in $(3.4)$ bounded by ${\epsilon_1 \over 2}x^v$. As a result, for each $x \in W_{ijp}$ there is an $\alpha$
with $0 \leq |\alpha| \leq m$ (or $2 \leq |\alpha| \leq M$) with
$$|x^{\alpha}\partial^{\alpha} f(x)| > \epsilon_1 f^*(x) \eqno (3.5)$$
Note that $(3.5)$ is independent of $i,j,$ or $p$ other than in the constant $\epsilon_1$. Thus by taking
$\epsilon_2$ to be the mininum of all such $\epsilon_1$, for any $x$ in some neighborhood of the origin 
there is an $\alpha$ with $0 \leq |\alpha| \leq m$ (or $2 \leq |\alpha| \leq M$) such that
$$|x^{\alpha}\partial^{\alpha} f(x)| > \epsilon_2 f^*(x) \eqno (3.6)$$
The above assumed $f(x)$ is real-analytic. In the case that $f(x)$ is merely smooth but with Newton 
polyhedron intersecting each coordinate axis,
the above argument gives $(3.6)$ if one replaces $f(x)$ by a finite Taylor polynomial $g(x)$ of $f(x)$ 
about the origin. But the Newton polyhedron condition ensures that if the expansion is taken out far
enough, the difference between $x^{\alpha}\partial^{\alpha} f(x)$ and $x^{\alpha}\partial^{\alpha} g(x)$
is less than $C |x| f^*(x)$ for any $\alpha$ appearing in $(3.6)$. Thus
by shrinking the neighborhood of the origin we are working in if necessary, once again $(3.6)$ holds and 
we are done with the proof of Lemma 3.6.

The next lemma shows each dyadic rectangle can be subdivided into boundedly many subrectangles such that
on each subrectangle, a single directional derivative of order $\leq m$ is bounded below by a constant
times the maximal possible value, so that one my apply the Van der Corput lemma in a single direction.

\noindent {\bf Lemma 3.7.} Suppose $f(x)$ is a smooth function satisfying $f(0) = 0$ and $\nabla f(0) = 0$
which has a nonvanishing
Taylor expansion at the origin, and that either $f(x)$ is real-analytic or a smooth function whose Newton polyhedron intersects
each coordinate axis. Again let $m$ denote the maximum 
order of any zero of any $f_F(x)$ and set $M =\max(2,m)$. 

Then there is a neighborhood $U$ of the origin and constants $K$ and $\eta$ such that if 
$R$ is any dyadic rectangle in $U$ then $R$ may divided into at most $K$ rectangles $R_j$ such that if
$2^{-k_i}$ denotes the length of $R$ in the $x_i$ direction, for each $R_j$ there is an $a$ and a single
$y =  (y_1,...,y_n)$ with  $|y_i| \leq 2^{-k_i}$ for all $i$ such that on $R_j$ we have
$$|(y \cdot \nabla)^a f(x)| \geq \eta |f^*(R_j)| \eqno (3.7)$$
Here $|f^*(R_j)|$ denotes $\sup_{R_j} f^*(x)$. For any $j$, in $(3.7)$ one can arrange that $0 \leq a 
\leq m$ or that $2 \leq a \leq M$.

\noindent {\bf Proof.} We consider only dyadic rectangles in the upper right octant as the analogous 
rectangles in other octants are treated the same way. Let $R$ be a dyadic rectangle $\prod_{i=1}^n [2^{-k_i + 1},
2^{-k_i}]$ contained in a neighborhood of 
the origin for which Lemma 3.6 is valid, and let $x \in R$. Then there is an $\alpha$ satisfying $0 \leq |\alpha| 
\leq m$ (or $2 \leq |\alpha| \leq M$) such that $(3.6)$ holds. As is well known (see [St] p.343 for
details), for a fixed $l$ there are a finite set of directions $\xi_1,...,\xi_{p_l} $ such that
every partial derivative of order $l$  can be written as a linear combination of the
$(\xi_i \cdot \nabla)^l$. Scaling this fact by a factor of $2^{k_i}$ in the $x_i$ variable, $(3.6)$ implies
that there is $\epsilon_3 > 0$ such that for any $x$ in a neighborhood of the origin there
exists some $y = (y_1,...,y_n)$ with $|y_i| \leq 2^{-k_i}$ for all $i$ such that at $x$ we have
$$|(y \cdot \nabla)^{|\alpha|} f(x)| > \epsilon_3 f^*(x)$$
Since $f^*(x)$ varies by a constant factor on $R$, we thus have
$$|(y \cdot \nabla)^{|\alpha|} f(x)| > \epsilon_4 |f^*(R)|\eqno (3.8)$$
The goal of this lemma then is to show that $R$ can be written as the union of boundedly
many subrectangles on each of which $(3.8)$ holds for a single $y$ and $|\alpha|$.
Suppose for some small $\delta > 0$, $z$ is such that $|z_i - x_i| \leq \delta2^{-k_i}$ for all $i$. Then
by the mean value theorem there is some $w$ between $x$ and $z$ such that
$$|(y \cdot \nabla)^{|\alpha|} f(z) - (y \cdot \nabla)^{|\alpha|} f(x)| \leq C\delta \sum_{|\beta| =
|\alpha| + 1}|w^{\beta}\partial^{\beta}f(w)| \eqno (3.9)$$
On the other hand, by Taylor expanding $x^{\beta}\partial^{\beta}f(x)$ about the origin, if $\sum_{\gamma} 
f_{\gamma}x^{\gamma}$ denotes the Taylor expansion of $f(x)$, then if $f(x)$ is real-analytic we have
$$|w^{\beta}\partial^{\beta}f(w)| \leq C\sum_{\gamma}|\gamma|^{|\alpha| + 1}|f_{\gamma}|w^{\gamma} 
\eqno (3.10)$$
By Theorem 3.2, this in turn is bounded by $C'|f^*(w)| \leq C''|f^*(R)|$ and we get
$$|w^{\beta}\partial^{\beta}f(w)| \leq C''|f^*(R)| \eqno (3.11)$$
In the case that $f(x)$ is smooth but with Newton polyhedron intersecting each axis, like in the previous
lemma $(3.11)$ holds for a Taylor polynomial approximation for $f(x)$ to high enough order and therefore 
for $f(x)$ as well. Combining $(3.9)$ and $(3.11)$, we get
$$|(y \cdot \nabla)^{|\alpha|} f(z) - (y \cdot \nabla)^{|\alpha|} f(x)| \leq C'''\delta|f^*(R)|  \eqno (3.12)$$
Thus if $\delta$ is smaller than ${\epsilon_4 \over 2C'''}$, for all $z$ such that $|z_i - x_i| \leq 
\delta2^{-k_i}$ for all $i$, $(3.8)$ and $(3.12)$ give
$$|(y \cdot \nabla)^{|\alpha|} f(z)| > {\epsilon_4 \over 2}|f^*(R)| \eqno (3.13)$$
Since $x$ was an arbitrary point in an arbitrary rectangle $R$, this completes the proof of Lemma 3.7. 

\noindent We are now in a position to prove Theorem 1.5. 

\noindent {\bf Proof of Theorem 1.5.} First, note that if $|\lambda_i|$ is larger than all other $|\lambda_j|$
for some $i < n+1$, in a small enough neighborhood of the origin one can integrate by parts in 
the $x_i$ variable and get that $|T(\lambda)| < C|\lambda_i|^{-1}$, which is better than what
we need here. So in what follows we always assume $|\lambda_{n+1}| \geq |\lambda_i|$ for all $i < n+1$.

We consider a small square centered at the origin on which Lemma 3.7 holds for $H(x)$, and assume $\Psi$ is
supported in this square. We divide the square
into dyadic rectangles and let $R = \prod_{i=1}^n [2^{-k_i + 1}, 2^{-k_i}]$ be one such rectangle.
Let $R = \cup_j R_j$ as in Lemma 3.7. Then the portion of $T(\lambda)$ coming from $R_j$ is given by
$$\int_{R_j} e^{-i\lambda_1 x_1-... - i\lambda_n x_n  - i\lambda_{n+1} H(x_1,...,x_n)} \Psi(x_1,...,x_n)
\,dx_1...\, dx_n \eqno (3.14)$$
We now scale the $x_i$ direction by $2^{k_i}$ to convert $R$ into $[{1 \over 2}, 1]^n$. Equation $(3.14)$
becomes
$$2^{-\sum_i k_i}\int_{R_j^*} e^{-i2^{-k_1}\lambda_1x_1...- i2^{-k_n}\lambda_n x_n  - i\lambda_{n+1}}
H(2^{-k_1}x_1,...,2^{-k_n}x_n)\Psi(2^{-k_1}x_1,...,2^{-k_n}x_n) \,dx \eqno (3.15)$$
Here $R_j^*$ is the scaled version of $R_j$. In view of Lemma 3.7 there is some $\delta > 0$ and
$y = (y_1,...,y_n)$ such that $|y_i| \leq 1$ for all $i$, such that for some $a$ with $2 \leq a \leq M$ 
($M = \max(m,2)$) on 
$R_j^*$ we have
$$|(y \cdot \nabla)^a H(2^{-k_1}x_1,...,2^{-k_n}x_n)| > \delta |H^*(R)| \eqno (3.16)$$
Hence we may apply the Van der Corput lemma in the $y$ direction and then integrate the result in the
$n-1$ orthogonal dimensions. The result is that $(3.14)$ is bounded by $C2^{-\sum_i k_i}\min(1, 
|\lambda_{n+1}|^{-{1 \over a}}|H^*(R)|^{-{1 \over a}})$. Adding this over all $j$, if $T_R(\lambda)$ denotes the portion
of $T(\lambda)$ coming from $R$, we have
$$|T_R(\lambda)| \leq C2^{-\sum_i k_i}\min(1, (|\lambda_{n+1}||H^*(R)|)^{-{1 \over a}})$$
$$\leq C 2^{-\sum_i k_i}\min(1, |\lambda_{n+1}|^{-{1 \over M}}|H^*(R)|^{-{1 \over M}}) \eqno (3.17)$$
Since $|H^*(R)|$ is within a constant factor of $H^*(x)$ for all $x$ and $2^{-\sum_i k_i}$ is a constant
factor times the area of $R$, we have that $(3.17)$ is bounded by
$$C \int_R \min(1, |\lambda_{n+1}|^{-{1 \over M}}|H^*(x)|^{-{1 \over M}}) \eqno (3.18)$$ 
Since we can assume $|\lambda_{n+1}| > |\lambda_i|$ for all $i < n+1$, we conclude that
$$|T_R(\lambda)| \leq C \int_R \min(1, |\lambda|^{-{1 \over M}}|H^*(x)|^{-{1 \over M}})\,dx \eqno (3.19)$$
Adding this over all $R$, if $S$ denotes the open square comprising the rectangles $R$ we have 
$$|T(\lambda)| \leq C \int_S\min(1, |\lambda|^{-{1 \over M}}|H^*(x)|^{-{1 \over M}}) \,dx$$
Shrinking $S$ if necessary, we assume $S$ is small enough that the decompositions of Theorems 3.1-3.3 hold
on $S$. Let $W_{ijp}$ be the sets coming from Theorems 3.1-3.3. Then the above is equal to
$$ C\sum_{ijp} \int_{W_{ijp}}\min(1, |\lambda|^{-{1 \over M}}|H^*(x)|^{-{1 \over M}}) \,dx\eqno (3.20)$$
In each term we can do the coordinate change $\beta_{ijp}:Z_{ijp} \rightarrow W_{ijp}$ given by Theorem 3.3.
By Theorem 3.3, each $x^v$ for vertex $v$ on the face $F_{ij}$ becomes the same monomial $m(z) = z^w$ in the 
$Z_{ijp}$ coordinates, and each $x^v$ for a vertex $v$ not on $F_{ij}$ becomes a monomial $z^{\alpha}$, 
$\alpha$ depending on $v$, such that $\alpha_l \geq w_l$ for all $l$. As a result we have 
$$H^* \circ \beta_{ijp}(z) \geq z^w$$
The Jacobian of $\beta_{ijp}(z)$ is a monomial, and we can compose this with a map of the form 
$a_{ijp}(z_1,...,z_n) = (z_1^{N_1},...,z_n^{N_n})$ so that  $\beta_{ijp} \circ a_{ijp}$ has constant
determinant. Letting $z^{w'}$ denote $m \circ a_{ijp}(z)$ and letting $Z_{ijp}'$ denote $\beta_{ijp}^{-1}
(Z_{ijp})$, a given term of $(3.20)$ is bounded by
$$C \int_{Z_{ijp}'}\min(1, (|\lambda| z^{w'})^{-{1 \over M}})\,dz \eqno (3.21)$$
By Theorem 2.6 of [G1], each component
$w_l'$ of $w'$ satisfies $0 \leq w_l' \leq d(H)$, where $d(H)$ denotes the Newton distance of $H$. The
theorem further says that at most $n - k$ of the $w_l'$ may be equal to $d(H)$, where $k$ denotes the
dimension of the face of $N(H)$ intersecting the critical line $\{(t,...,t): t > 0\}$ in its interior.

We now can directly apply Lemma 3.4 to bound $(3.21)$ and complete the proof of Theorem 1.5. 
We start with part a). In this case, $d(H) < 2$ and $M = 2$. Since $d(H) < 2$, each $w_i' < 2$. Consequently, 
${w_l' \over M} < 1$ for all $l$. As a result, parts a) and b) of Lemma 3.4 apply. Part b) gives the larger 
term,
and we conclude that each term of $(3.21)$ is bounded by $C|\lambda|^{-{1 \over 2}}$. Adding over all $i,j$,
and $p$ gives that $|T(\lambda)| < C|\lambda|^{-{1 \over 2}}$. This gives part a) of Theorem 1.5.

Next, we move to part b). In this case $M \leq d(H)$ and $d(H) > 2$. Actually, since $M$ is an integer we 
have $M \leq \lfloor d(H) \rfloor$. Thus $(3.21)$ is at most
$$ C \int_{Z_{ijp}'}\min(1, (|\lambda| z^{w'})^{-{1 \over \lfloor d(H) \rfloor}})\,dz \eqno (3.22)$$
First we consider the case where $d(H)$ is an integer. Since by above $0 \leq w_l' \leq d(H)$ with equality 
holding for at most $n - k$ indices, we have that ${w_l' \over \lfloor d(H) \rfloor}  = {w_l' \over  d(H)} \leq 1$ for 
all $l$, with equality holding for at most $n - k$ values of $l$. We now add parts a) and c) of Lemma 3.4,
and add the result over all $i,j$, and $p$. We get 
$$|T(\lambda)| \leq C'|\lambda|^{-{1 \over d(H)}}(\ln |\lambda|)^{n - k} \eqno (3.23)$$
This is the estimate we seek. Next, we suppose $d(H)$ is not an integer and thus $\lfloor d(H) \rfloor
< d(H)$. Now, for a given $(i,j,p)$, it is possible that some ${w_l' \over \lfloor d(H) \rfloor}
> 1$. In this case, we add parts a) and d) of Lemma 3.4 and get that 
$(3.22)$ is at most $C'|\lambda|^{-{1 \over  \lfloor d(H) \rfloor }}(\ln |\lambda|)^{n- k - 1}$, which is
less than $C'|\lambda|^{-{1 \over  d(H)}}(\ln |\lambda|)^{n- k - 1}$, the desired upper bounds.
If on the other hand each ${w_l' \over \lfloor d(H) \rfloor} \leq 1$, we add part a) of Lemma 3.4 
to either part b) or c) of that theorem, depending on whether or not the maximal 
${w_l' \over \lfloor d(H) \rfloor}$ is less than or equal to one. In either case, one  gets that
$(3.22)$ is at most  $C'|\lambda|^{-{1 \over \lfloor d(H) \rfloor }}(\ln |\lambda|)^{n- k}$
which is bounded by the desired upper bound $C''|\lambda|^{-{1 \over  d(H)}}(\ln |\lambda|)
^{n- k - 1}$. Adding over all $i$, $j$, and $p$ gives the desired upper bounds for Theorem 1.5 b) and
we have proven part b). 

We now move to part c). Each $w_l' \leq d(H)$ and $d(H) < M = m$, so 
like in part a), each exponent ${w_l' \over M}$ is less than 1. Hence we again add parts a) and b) of 
Lemma 3.4, with part b) giving the larger term, and adding over all $i,j,$ and $p$ we get
$$|T(\lambda)| \leq C'|\lambda|^{-{1 \over M}} \eqno (3.24)$$
This gives part c) of Theorem 1.5 and we are done.

\noindent Lastly, we give the proof of Theorem 1.6. 

\noindent {\bf Proof of Theorem 1.6.} 

We start with part a). By Theorem 1.6a) of [G1], whenever $\phi$ 
is nonnegative with $\phi(y) \neq 0$, $\limsup_{\lambda_{n+1} \rightarrow \infty} 
{T((0,...,0,\lambda_{n+1})) \over |\lambda_{n+1}|^{1 \over d(H)}\ln|\lambda_{n+1}|^{n-k-1}} > 0$. This
immediately gives part a). As for b), suppose $V$ is any neighborhood of the origin and suppose $a^k \in V$. Then
after doing a coordinate change from turning $x_i - b_k(x_1,...x_{i-1},x_{i+1},...,x_n)$ into $x_i$ and 
leaving the other variables fixed, we may assume $H(x)$ has the phase $\xi_k(x)x_i^m$. In order to prove
Theorem 1.6b) it suffices to show that in the new coordinates, on any sufficiently small neighborhood $W$ of 
$\bar{a}^k = (a^k_1,...,a^k_{i-1},0,a^k_{i+1},...,a^k_n)$ there is a bump function $b(x)$ supported in $W$ such 
that $\bar{T}(\lambda)$ satisfies the conclusions of Theorem 1.6b), where 
$$\bar{T}(\lambda) = \int e^{-i\lambda_1 x_1... - i\lambda_n x_n  - i\lambda_{n+1} \xi_k(x)x_i^m} 
b(x_1,...,x_n)\,dx_1...\,dx_n \eqno (3.25)$$
Since $\xi_k(\bar{a}^k) \neq 0$ in the new coordinates, after another coordinate change we can assume 
$\xi_k(x)x_i^m$
is actually $x_i^m$ and $b(x_1,...,x_n)$ is a different cutoff function but one we may still choose 
freely. In particular, we may let $b(x_1,...,x_n) = b_1(x_i)b_2(x_1,...,x_{i-1},x_{i+1},...,x_n)$, 
where $\hat{b}_2(0) \neq 0$ and where $\int b_1(x)e^{-i\lambda x^m}\,dx$ decays as $C_0|\lambda|^{-{1 \over 
m}}$ as $|\lambda| \rightarrow \infty$, where $C_0$ is nonzero. Then $(3.25)$ becomes
$$\bar{T}(\lambda) = \hat{b}_2(\lambda_1,...,\lambda_{i-1},\lambda_{i+1},...,\lambda_n) \int
b_1(x)e^{-i\lambda_i x^m}\,dx$$
Setting $\lambda = (0,..,0,\lambda_i,0,...,0)$ gives decay of $C_0|\lambda|^{-{1 \over 
m}}$ as $|\lambda| \rightarrow \infty$, giving part b) and we are done. 

\noindent {\bf 4. Maximal operators: some lemmas and the proof of the $L^\infty$ bounds.}

We start by showing that when the Hessian condition of Theorem 1.1 or 1.3 is not satisfied, one
has stringent restrictions on what the $H_F(x)$ may be. This ensures that in examples 1 and 2
of section 1, the Hessian condition is automatically satisfied, and that in Theorem 1.1b) it is satisfied
unless $N(H)$ has the special form given there.

\noindent {\bf Lemma 4.1.} Suppose $f(x)$ is a smooth function satisfying $f(0) = 0$ and $\nabla f(0) = 0$
which has a nonvanishing
Taylor expansion at the origin. Further suppose that there do not exist directions
$u$ and $v$ such that the Hessian of $f(x)$ in the $u$ and $v$ variables vanishes to infinite order at the
origin. Then there is an $m > 1$ such that each vertex $v$ of $N(f)$ lies on a coordinate axis at 
height $m$. Furthermore, if $F$ is any compact face of $N(f)$, $f_F(x)$ is of the form 
$c_F(\sum_{i=1}^n \beta_i^F x_i)^m$ where $c_F \neq 0$, but some $\beta_i^F$ may be zero. The number $m$ is greater
than the Newton distance $d(f)$ unless the vertex set of $N(f)$ consists of a single vertex lying on a 
coordinate axis. Hence in case b) Theorem 1.1, the Hessian condition is automatically satisfied whenever
the vertex set of $N(f)$ does not solely consist of a single vertex lying on a coordinate axis.

\noindent {\bf Proof.} Let $F$ be any compact face of $N(f)$. Write the 
Taylor expansion of $f(x)$ at the origin as $\sum_{\alpha}f_{\alpha}x^{\alpha}$  $= f_F(x) + \sum_{\alpha 
\notin F}f_{\alpha}x^{\alpha}$.   
Since the Hessian condition is assumed to hold for no $u$ and $v$, it in particular does not hold if $u$
and $v$ are coordinate directions $x_i$ and $x_j$. Letting $H_{ij}f(x)$ and $H_{ij}f_F(x)$ denote the 
Hessian matrices of $f$ and $f_F$ respectively in 
the $x_i$ and $x_j$ variables, write $H_{ij}f(x) = H_{ij}f_F(x) + E(x)$.  
Let $H_{ij}^*f(x)$ denote the matrix obtained from $H_{ij}f(x)$ by multiplying the $ii$ entry by $x_i^2$,
the $jj$ entry by $x_j^2$, and the $ij$ and $ji$ entries by $x_ix_j$, with analogous definitions for
$H_{ij}^*f_F(x)$ and $E^*(x)$. Then we have
$$H_{ij}^*f(x) = H_{ij}^*f_F(x) + E^*(x) \eqno (4.1)$$
We may let $c = (c_1,...,c_n)$ be a vector with positive entries such that for $w = (w_1,...,w_n) \in N(f)$, 
$c \cdot w$ is minimized exactly when $w \in F$. Denote this minimal value of $c \cdot w$ by $a$.
Note that any term $m_{\alpha}x^{\alpha}$ appearing in $H_{ij}^*f_F(x)$ satisfies $c \cdot \alpha = a$,
while each term $m_{\alpha}x^{\alpha}$ appearing in the Taylor expansion of an entry of $E^*(x)$ satisfies  $c \cdot \alpha > a$. As a 
result, in view of $(4.1)$ the determinant of $H_{ij}^*f(x)$ consists of the determinant of $H_{ij}^*f_F(x)$,
a polynomial for which any nonzero term $m_{\alpha}x^{\alpha}$ satisfies $c \cdot \alpha = 2a$, plus a smooth
function for which any term $m_{\alpha}x^{\alpha}$ of its Taylor expansion satisfies $c \cdot \alpha > 2a$.
Since the determinant of $H_{ij}f(x)$ is assumed to vanish to infinite order at the origin, the same is
true of the determinant of $H_{ij}^*f(x)$. Hence the terms of its Taylor expansion satifiying $c \cdot 
\alpha = 2a$ must all be zero. Equivalently, the determinant of $H_{ij}^*f_F(x)$ is identically zero, which
means the same is true for the determinant of $H_{ij}f_F(x)$.

The indices $i$ and $j$ were arbitrarily chosen, so we conclude that for all $i$ and $j$ the determinant of $H_{ij}f_F(x)$
is the zero function. This can only happen if the $n$ by $n$ Hessian matrix of $f_F(x)$ has rank $\leq 1$
everywhere. Since $f_F(x)$ is a polynomial, as described after the statement of Theorem 1.1 this means 
that $f_F(x)$ is of the form $L(x) + p(\sum_i\beta_i x_i)$,
where $L$ is linear and $p$ is a polynomial of degree at least 2. Since $f(x)$ has a critical point at the origin, so does
$f_F(x)$ and thus $L = 0$. Since each term $f_{\alpha}x^{\alpha}$ of $f_F(x)$'s Taylor expansion 
satisfies $c \cdot \alpha = a$, $p(t)$ must be a monomial; otherwise $c \cdot \alpha$ would take on 
various multiples of its minimal value. Hence we conclude that $f_F(x)$ is of the desired form 
$c_F(\sum_{i=1}^n\beta_i^F x_i)^{m_F}$.

Next, we show that any vertex $v$ of $N(f)$ lies on one of the coordinate axes. To see this, observe that
since $v$ is a face of $N(f)$, the above applies to $v$ and we have that the monomial $f_vx^v$ must be of 
the form $c_v(\sum_{i=1}^n\beta_i^v x_i)^{m_v}$. Clearly this can only happen if there is exactly one 
nonvanishing $\beta_i$, in which case $v$ lies on the $x_i$ axis. 

To finish the proof of the lemma, we will show that if $N(f)$ has more than one vertex, then each $m_F$
is a single number $m$ such that $m > d(f)$. Suppose $v_1$ and $v_2$ are distinct vertices of $N(f)$. 
The segment $l$ 
connecting $v_1$ to $v_2$ is also a face of $N(f)$, and so $f_l(x)$ is also of the form $c_l(\sum_
{i=1}^n\beta_i^l x_i)^{m_l}$. Since it connects two vertices lying on two separate coordinate axes, 
$f_l(x)$ can be written as  
$$f_l(x) = c_l(\beta_i^lx_i + \beta_j^lx_j)^{m_l} \eqno (4.2)$$
Here $c_l$, $\beta_i^l$, and $\beta_j^l$ are all nonzero. Note that $(4.2)$ forces $v_1$ and $v_2$ to be
at the same height $m_l$ on their respective axes. Since $v_1$ and $v_2$ were arbitrary, this means that
there is a single $m$ such that $m = m_v$ for all $v \in v(f)$. This in turn implies that $m = m_F$ for all
$F$ since each vertex of $F$ will be at height $m$ on its coordinate axis. Since $f$ and its gradient are
both zero at the origin, $m > 1$. Lastly, we note that the
Newton polyhedron $P$ generated by the two vertices $v_1$ and $v_2$ alone satisfies $d(P) = {m \over 2}$, so 
$d(f)$ is at most this value. Hence $d(f) \leq {m \over 2} < m$ and we are done. 

\noindent {\bf Lemma 4.2.} Suppose $p(x)$ and $q(x)$ are smooth functions satisfying $p(0) = 0$, $q(0) = 0$, 
$\nabla p(0) = 0$, and $\nabla q(0) = 0$. Assume also that $p$ and $q$ have nonvanishing Taylor expansion at the 
origin.  Suppose further that there is $r > 0$ such that
$N(q) = rN(p)$; that is, $\alpha \in N(p)$ iff $r\alpha \in N(q)$. Then $N(pq) = (r + 1)N(p)$, and 
furthermore if $F$ is a compact face of $N(p)$ then $pq_{(r+1)F}(x) = p_F(x)q_{rF}(x)$. Here $rF$ denotes the
dilation of the face $F$ by $r$ and $(r+1)F$ the dilation by $r + 1$. 

\noindent {\bf Proof.} Let $\sum_{\alpha} p_{\alpha}x^{\alpha}$ denote the Taylor expansion of $p(x)$ at the
origin. Let $c = (c_1,...,c_n)$ be a vector with positive entries. Then there is some face $F$ of $N(p)$ 
such that $c \cdot \alpha$ achieves its minimal value amongst $\alpha$ with $p_{\alpha} \neq 0$ iff 
$\alpha \in F$.
We write $p(x) = p_F(x) + E_p(x)$, and similarly write $q(x) = q_{rF}(x) + E_q(x)$. 
The terms of $p_F(x)$ and $q_{rF}(x)$ minimize $c \cdot \alpha$ in the Taylor expansions of $p(x)$ and 
$q(x)$ respectively. If for $p_F(x)$ they satisfy
$c \cdot \alpha = a$, then for $q_{rF}(x)$ they satisfy $c \cdot \alpha = ar$. Hence multiplying together, one
gets $p(x) = p_F(x)q_{rF}(x) + E_{pq}(x)$, where the terms of $p_F(x)q_{rF}(x)$ now satisfy $c \cdot \alpha = 
(a + 1)r$, and each monomial appearing in  $E_{pq}(x)$ satisfies $c \cdot \alpha > (a + 1)r$. Hence the 
terms in the Taylor expansion of $p(x)q(x)$ minimizing $c \cdot \alpha$ are those of $p_F(x)q_{rF}(x)$, 
and $(r+1)F$ is the corresponding face of $N(pq)$ with $pq_{(r+1)F}(x) = p_F(x)q_{rF}(x)$. Since $c$ was 
arbitrary, we are done. 

Recall that in section 2, we defined functions $H^*(x) = \sum_{v \in v(H)} |x|^v$ and 
$H^{**}(x) = \big(\sum_{i,j} x_i^2 x_j^2 (\partial_{x_ix_j}^2 H(x))^2\big)^{1 \over 4}$ for use in our 
damping functions. The next lemma records some relevant properties of the latter function.

\noindent {\bf Lemma 4.3.} Let $\bar{H}(x) = H^{**}(x)^4 = \sum_{i,j} x_i^2 x_j^2 (\partial_{x_ix_j}^2 
H(x))^2$. Then $N(\bar{H}) = 2N(H)$. Let $m$ denote the maximum order of any zero of any $H_F(x)$ in 
$(\R - \{0\})^n$. Then the maximum order of any zero of any $\bar{H}_F(x)$ in $(\R - \{0\})^n$ is 
$2\max(m,2) - 4$. 

\noindent {\bf Proof.} If $\sum_{\alpha} H_{\alpha}x^{\alpha}$ denotes the Taylor 
expansion of $H(x)$ at the origin, then if $i \neq j$, $x_ix_j\partial_{x_ix_j}^2 H(x)$ has Taylor expansion 
$\sum_{\alpha} \alpha_i\alpha_jH_{\alpha}x^{\alpha}$, and if $i = j$ then $x_ix_j\partial_{x_ix_j}^2 H(x)$
has Taylor expansion $\sum_{\alpha} \alpha_i(\alpha_i - 1)H_{\alpha}x^{\alpha}$. In both cases, 
$N(x_ix_j\partial_{x_ix_j}^2 H) \subset N(H)$. Hence by Lemma 4.2, $N((x_ix_j\partial_{x_ix_j}^2 H)^2) 
\subset 2N(H)$. Adding over all $i$ and $j$, we get $N(\bar{H}) \subset 2N(H)$. 

To show equality, let $v$
be a vertex of $N(H)$ and let $H_vx^v$ denote the corresponding term in $H$'s Taylor expansion.
Since $H(x)$ has a zero of order at least 2 at the origin there is either an $i$ such that $v_i \geq 2$ and therefore
$v_i(v_i - 1) \neq 0$, or there are $i$ and $j$ such that $v_i, v_j \geq 1$ and therefore $v_iv_j \neq 0$. 
Thus $v$ will show up in the Taylor expansion of at least one $x_ix_j\partial_{x_ix_j}^2 H$. Using 
Lemma 4.2 again, this means
$2v$ will show up in in the Taylor expansion of $x_i^2 x_j^2(\partial_{x_ix_j}^2 H)^2$, and with positive
coefficient. Adding over all $i$ and $j$ and using that Lemma 4.2 implies that $2v$ can never appear with 
negative coefficient, we conclude $2v$ 
will appear in the Taylor expansion of $\bar{H}(x)$ and therefore  $2v \in N(\bar{H})$. Taking the convex hull over
all $v_i$ gives that $2N(H) \subset N(\bar{H})$. We conclude that $2N(H) = N(\bar{H})$ as needed.

We now move to the statement about the orders of the zeroes of the $\bar{H}_F(x)$. Suppose $x \in 
(\R - \{0\})^n$ is such that some $H_F(x)$ is nonzero or has a zero of order $\leq 2$ 
at $x$. Then by Lemma 3.5 there are some $i$ and $j$ such that $\partial_{x_ix_j}^2 H_F(x) \neq 0$. Hence 
$x_ix_j\partial_{x_ix_j}^2 H_F(x) \neq 0$. By Lemma 4.2, $\big((x_ix_j\partial_{x_ix_j}^2 H)^2\big)_{2F}(x) =
[(x_ix_j\partial_{x_ix_j}^2 H)_F(x)]^2$. This in turn is equal to  $[x_ix_j\partial_{x_ix_j}^2 H_F(x)]^2$, 
a positive quantity. Adding over all $i$ and $j$, we get $\bar{H}_
{2F}(x) > 0$. Hence $\bar{H}_{2F}$ has a zero of order $0 \leq 2\max(m,2) - 4$ at $x$. In the event that
$m \leq 2$, we can conclude that each $\bar{H}_{2F}$ is nonvanishing on $(\R - \{0\})^n$, completing the 
proof of the lemma for when $m \leq 2$.

On the other hand, suppose $x \in (\R - \{0\})^n$ and $F$ are such that $H_F$ has a zero of order $o > 2$
at $x$. Note that necessarily $o \leq m$. Then each $\partial_{x_ix_j}^2 H_F$ has a zero of order at least
$o - 2$ at $x$, with at least one having order exactly $o - 2$. Hence the same is true for each 
$x_ix_j\partial_{x_ix_j}^2 H_F$. Adding up the squares of these functions as in the previous paragraph,
we obtain that $\bar{H}_{2F}$ has
a zero of order $2o - 4 \leq 2m - 4 = 2\max(m,2) - 4$ at $x$. If $x$ is such that $o = m$, then we have
equality. This completes the proof of Lemma 4.3.

\noindent We are now in a position to set up the interpolation for Theorems 1.1 and 1.3. 

\noindent {\bf The interpolation of Theorems 1.1 and 1.3.} 

In view of Lemma 4.1, under the hypotheses of Theorem 1.1 or 1.3 we can assume there exist directions $u$ and $v$ such that the
determinant $D(x)$ of the Hessian of $H$ in the $u$ and $v$ directions does not vanish to infinite order at
the origin. As above, we let $H^*(x_1,...,x_n) = \sum_{v \in v(H)} |x|^
v$ and $H^{**}(x) = \big(\sum_{i,j} x_i^2 x_j^2 (\partial_{x_ix_j}^2 H(x))^2\big)^{1 \over 4}$. Then 
for a small $\delta > 0$, if $d$ denotes $d(H)$ we define the 
damping function $P(x)$ by
$$P(x) = |D(x)|^{\delta}H^*(x)^{-{1 \over \max(d,2)}}H^{**}(x) \eqno (4.3)$$
We then define the damped surface measure $d\sigma_z(s) = |P(x)|^z d\sigma(s)$ and denote the 
maximal operator associated to $e^{z^2}\sigma_z$ by $M_z$. Our next lemma gives the needed $L^{\infty}$ 
boundedness properties of the $M_z$.

\noindent {\bf Lemma 4.4.} Let $m$ denote the maximum order of any zero of any $H_F(x)$, and let $M = 
\max(2, m)$. Suppose one of the following holds.

\noindent (a) $d > 2$ and $Re(z) > - \min({2 \over M - 2}, {2 \over d - 2})$. (Taken as $-{2 \over d - 2}$
if $M = 2$). \hfill \break
\noindent (b) $d \leq 2$, $M > 2$, and $Re(z) > -{2 \over M - 2}$. \hfill \break
\noindent (c) $d \leq 2, M = 2$ and $z$ is arbitrary.

\noindent Then if $\delta$ in 
$(4.3)$ is suffciently small, there is a constant $A$ depending on $Re(z)$ such that $|\sigma_z|(S) \leq
A$. Consequently, as an operator on $L^{\infty}(\R^n)$, $||M_z|| \leq A$.

\noindent {\bf Proof.} We apply Lemma 3.7 to the function $\bar{H}(x)$. Correspondingly, we divide some cube
$[-2^{-K_1},2^{-K_1}]^n$ into dyadic rectangles $R$, each of which we further subdivide into boundedly many rectangles 
$R_j$ on which Lemma 3.7 holds. Define $\bar{H}^*(x) = \sum_{v \in N(\bar{H})} |x|^v$ and
$\bar{H}^*(R) = 
\sup_{R} \bar{H}^*(x)$. By Lemma 3.7, for some $\delta > 0$ there is an $a$ with $0 \leq a \leq 2M  - 4$ 
and a $y = (y_1,...,y_n)$ with $|y_i| \leq 2^{-k_i}$ for all $i$ such that on $R_j$ we have
$$|(y \cdot \nabla)^a \bar{H}(x)| \geq \delta \bar{H}^*(R)\eqno (4.4)$$
Note that since $N(H^*) = 2N(H)$, $\bar{H}^*(R)$ is comparable to $(H^*(R))^2$. Hence $(4.4)$ implies
$$|(y \cdot \nabla)^a \bar{H}(x)| \geq \delta' (H^*(R))^2 \eqno (4.4')$$
We first consider the case where $a = 0$. Define  $P_0(x) = (H^*(x))^{-{1 \over \max(d,2)}} H^{**}(x)$. 
Since $H^{**}(x) = \bar{H}(x)^{1 \over 4}$, 
$(4.4')$ implies that on $R_j$, for $t > 0$ we have 
$$P_0(x)^{-t} \leq C|H^*(R)|^{-t({1 \over 2} - {1 \over \max(2,d)})}$$
Integrating this over $R_j$ gives
$$\int_{R_j}P_0(x)^{-t} \leq C|R_j| |H^*(R)|^{-t({1 \over 2} - {1 \over \max(2,d)})} \eqno (4.5)$$
This is the estimate we will need for when $a = 0$. Now suppose $a > 0$. We will again bound 
$\int_{R_j}P_0(x)^{-t}$. Without loss of generality we assume 
$R$ is in the upper right octant and we write $R = \prod_{i=1}^n [2^{-k_i - 1}, 2^{-k_i}]$. We change 
coordinates in the integral $\int_{R_j}P_0(x)^{-t}$, scaling the $x_i$ coordinate by $2^{k_i}$ so that
$R$ becomes the cube $\prod_{i=1}^n [{1 \over 2}, 1]$. Letting $R_j^*$ denote $R_j$ in the new coordinates, we
have
$$\int_{R_j}P_0(x)^{-t} = 2^{-\sum_i k_i}\int_{R_j^*}P_0(2^{-k_1}x_1,...,2^{-k_n}x_n)^{-t} \eqno (4.6)$$
Since $|H^*(x)|$ is within a constant of the fixed value $H^*(R)$ on $R_j$, we have
$$\int_{R_j}P_0(x)^{-t} \leq C2^{-\sum_i k_i}|H^*(R)|^{t \over \max(2,d)} \int_{R_j^*}|H^{**}(2^{-k_1}x_1
,...,2^{-k_n}x_n)|^{-t}$$
$$ = C2^{-\sum_i k_i}|H^*(R)|^{t \over \max(2,d)} \int_{R_j^*}|\bar{H}(2^{-k_1}x_1
,...,2^{-k_n}x_n)|^{-{t \over 4}} \eqno (4.7)$$
(We of course don't need absolute values in $(4.7)$ since the quantity is nonnegative, but we include them
anyhow for readability).
The condition $(4.4')$, translated into the new coordinates, says that for some $y = (y_1,...,y_n)$ with
$|y_i| \leq 1$ for all $i$ we have 
$$|(y \cdot \nabla)^a [\bar{H}(2^{-k_1}x_1,...,2^{-k_n}x_n)]| > \delta' (H^*(R))^2 \eqno (4.8)$$
Hence by the one-dimensional Van der Corput theorem for measures (see [Ch]), for any $\epsilon > 0$ we have 
$$|\{x \in R_j^*: |\bar{H}(2^{-k_1}x_1,...,2^{-k_n}x_n)| < \epsilon\}| \leq C 
|\{x \in [0,1]^n: \delta' (H^*(R))^2 x_1^a< \epsilon\}| \eqno (4.9)$$
As a result, by the relation between $L^p$ norms and distribution functions, we have 
$$\int_{R_j^*}|\bar{H}(2^{-k_1}x_1,...,2^{-k_n}x_n)|^{-{t \over 4}} \leq C
\int_{[0,1]^n} (\delta' (H^*(R))^2 x_1^a)^{-{t \over 4}} \eqno (4.10)$$
The right-hand side is finite whenever $t < {4 \over a}$. In particular, by Lemma 4.3 it is finite 
whenever $t < {2 \over M - 2}$. In this case, $(4.10)$ is bounded by $C_0 (H^*(R))^
{-{t \over 2}}$, where $C_0$ depends on $t$ as well as the function $H(x)$. Putting this back into
$(4.7)$, we get that 
$$\int_{R_j}P_0(x)^{-t} \leq C_12^{-\sum_i k_i}|H^*(R)|^{-t({1 \over 2} - {1 \over \max(2,d)})}$$
$$\leq C_1 |R_j||H^*(R)|^{-t({1 \over 2} - {1 \over \max(2,d)})} \eqno (4.11)$$
Here $C_1$ depends on $t$ and the function $H(x)$. Note that $(4.11)$ and $(4.5)$ are the same other than
the constants. Hence we may add $(4.5)$ or $(4.11)$ over all $R_j$ comprising a given $R$,
and we get
$$\int_R P_0(x)^{-t} \leq C_2 |R||H^*(R)|^{-t({1 \over 2} - {1 \over \max(2,d)})}$$
$$\leq C_3 \int_R |H^*(x)|^{-t({1 \over 2} - {1 \over \max(2,d)})} \eqno (4.12)$$
The last inequality follows from the fact that $H(x)$ is within a constant of a fixed value on $R$, and
once again $C_3$ depends on $t$ and the function $H(x)$. We now add $(4.12)$ over all $R$ to obtain
$$\int P_0(x)^{-t} \leq C_3 \int |H^*(x)|^{-t({1 \over 2} - {1 \over \max(2,d)})} \eqno (4.13)$$
If $d \leq 2$, the integral $(4.12)$ is always
finite. If $d > 2$, Theorem 1.2 of [G1] says that $|H^*(x)|^{-u}$ is integrable near the origin iff $u$ is 
less than ${1 \over d(H^*)} = {1 \over d}$. Hence $(4.12)$ is finite whenever
$t < {{1 \over d}  \over {1 \over 2} - {1 \over d}} = {2 \over d - 2}$. The other condition we had
for finiteness of  $\int_R P_0(x)^{-t}$ came from $(4.10)$, which we 
saw was finite whenever $t < {2 \over M - 2}$. This restriction only arose when $a$ could be positive,
which can only happen if $M > 2$. Hence if $d > 2$, 
$\int_R P_0(x)^{-t}$ is finite whenever $t < \min({2 \over M - 2}, {2 \over d - 2})$, while if $d \leq 2$, we 
have one restriction $t < {2 \over M - 2}$, occurring when $M > 2$. Note that these $t$ are the exponents of this
lemma. So for any such $t$, by Holder's inequality if $\delta > 0$ is small enough, we have
$$\int (P_0(x)|D(x)|^{\delta})^{-t} < C_4 \eqno (4.14)$$
Since $P(x) = P_0(x)|D(x)|^{\delta}$ and thus $|P(x)^z| = (P_0(x)|D(x)|^{\delta})^{Re(z)}$, we 
conclude that the measure $|\sigma_z|$ is uniformly bounded in $|Im(z)|$ whenever $Re(z)$ satisfies the 
conditions of the lemma. This completes the proof of Lemma 4.4. 

\noindent {\bf 5. Maximal Operators: $L^2$ bounds and the proofs of Theorems 1.1 and 1.3.}

The main goal of this section are the following estimates on the Fourier transform of the measures 
$\sigma_z$ that will enable us to invoke the lemma of Sogge-Stein.

\noindent {\bf Theorem 5.1}. If $\delta$ is chosen sufficiently small and $Re(z) > 1$, there are 
constants $A$ and $\epsilon$ depending on $Re(z)$, $H(x)$, and $\delta$ such that 
$$|\hat{\sigma}_z(\lambda)| < A(1 + |Im(z)|)(1 + |\lambda|)^{-{1 \over 2} - \epsilon} \eqno (5.1a)$$
In addition, for all $i$ we have
$$|\partial_{\lambda_i} \hat{\sigma}_z(\lambda)| < A(1 + |Im(z)|)(1 + |\lambda|)^{-{1 \over 2} - \epsilon}
\eqno (5.1b)$$

\noindent {\bf Proof.} Our focus will be on proving $(5.1a)$ as the proof of $(5.1b)$ is identical with
$\Psi(x)$ replaced by $x_i\Psi(x)$ or $H(x)\Psi(x)$. As explained in section 2 above $(2.3)$, our task is 
to show that if $\delta$ is sufficiently small and $Re(z) > 1$, then 
$|G_z(\lambda)|$ is bounded by the right-hand side of $(5.1a)$, where
$$G_z(\lambda) = \int e^{-i\lambda_1 x_1... - i\lambda_n x_n  - i\lambda_{n+1} H(x)} |D(x)|^{\delta z}
|H^*(x)|^{-{z \over \max(d,2)}} |H^{**}(x)|^z\Psi(x)\,dx_1...\,dx_n \eqno (5.2)$$ 
(Once again, although we don't need absolute values on $H^*(x)$ or $H^{**}(x)$ since they are nonnegative,
we include them in our arguments to improve readability).

We divide the domain of integration into dyadic rectangles, and without loss of generality we consider
only those in the upper right octant as the other octants are dealt with in exactly the same way. We further
may consider only those dyadic rectangles whose shortest side is at least $|\lambda|^{-1}$, as the measure
of the union of the remaining rectangles is at most $C|\lambda|^{-1}$ and therefore will not affect the
truth of $(5.1a)$. Since there are $C(\log |\lambda|)^n$ rectangles remaining, it suffices to prove
that the portion of $G_z(\lambda)$ over a given rectangle is bounded by the right-hand side of $(5.1a)$.
Similarly, if $N_0$ is such that $|\{x: |H^*(x)| < |\lambda|^{-N_0}\}| < C|\lambda|^{-1}$, it suffices
to consider only those rectangles $R$ such that $\sup_R |H^*(x)| > |\lambda|^{-N_0}$. 

Next for sufficiently small positive $c$ and $\delta_1$, we divide each remaining dyadic rectangle $R$ into 
subrectangles of radius $c|\lambda|^{-{\delta_1}}$. We will show that the portion of $G_z(\lambda)$ coming
from each such subrectangle is bounded by the right-hand side of $(5.1a)$, and $\epsilon$ will not 
depend on $\delta_1$. Furthermore, the constant $c$ will be independent of $Im(z)$. Hence letting 
$\delta_1$ be small enough and adding over the subrectangles comprising $R$, we will have $(5.1a)$.
Thus we let $S$ be any one of these subrectangles and focus our attention on bounding the part of the 
integral over $S$. 

We do this transfer from the rectangles $R$ to the subrectangles $S$ so that we may effectively replace $H(x)$, 
$D(x)$, and $H^{**}(x)$ in $(5.2)$ by polynomial approximations of sufficiently high degree. This will 
enable us to apply the proof of the Van der Corput lemma which requires the integrals involved to be over 
boundedly many intervals on each of which a relevant differentiated function is monotone.
(If these functions were all polynomials to begin with, we would not have to do this transfer). To this end,
we fix $x_0 \in S$ and for a given $N$ let $H_N(x)$ be the finite Taylor approximation to $H(x)$ about $x_0$.
Let $G_z^S(\lambda)$ be the portion of $(5.2)$ over $S$, and let $(G_z^S)'(\lambda)$ be this integral
with $H(x)$ replaced by $H_N(x)$. Then $|G_z^S(\lambda) - (G_z^S)'(\lambda)|$ is bounded by
$$ \int |e^{-i\lambda_{n+1}(H(x) - H_N(x))} - 1|
|D(x)|^{\delta Re(z)}|H^*(x)|^{-{1 \over \max(d,2)}Re(z)} |H^{**}(x)|^{Re(z)}|\Psi(x)|\,dx_1...\,dx_n \eqno (5.3)$$
In $(5.8a)$ below, we will see that $|H^*(x)|^{-{1 \over \max(d,2)}} |H^{**}(x)|$ is bounded. 
So if $N > {2 \over \delta_1}$ and $c$ is sufficiently small, $(5.3)$ shows that $|G_z^S(\lambda) - (G_z^S)'
(\lambda)| \leq C|\lambda|^{-1}$, better than the estimate we need. So we may replace $H(x)$ by $H_N(x)$ in 
our future arguments. 

Next, we do a polynomial approximation to the $H^{**}(x)$ appearing in $(5.3)$. (We do
not have to do anything with the $|H^*(x)|$ factor since it is already a polynomial when restricted to $S$).
Namely, we define $H^{**}_N(x)$ by
$$H^{**}_N(x)= \big(\sum_{i,j} x_i^2 x_j^2 (\partial_{x_ix_j}^2 H_N(x))^2\big)^{1 \over 4}$$
Let
$\bar{H}(x) = H^{**}(x)^4$ and $\bar{H}_N(x) = H^{**}_N(x)^4$. Then $\bar{H}(x) - \bar{H}_N(x)$ has a zero 
of order at least $2N$ at $x_0$. The difference between $(G_z^S)'(\lambda)$ and the expression one gets if one
replaces $H^{**}(x)$ by $H^{**}_N(x)$ in the integral is bounded by
$$C\int_S |H^*(x)|^{-{1 \over \max(d,2)}Re(z)}|\bar{H}(x)^{z \over 4} - \bar{H}_N(x)^{z \over 4}| \eqno (5.4)$$
Recall there is some $N_0$ such that we are assuming $|H^*(x)| \geq |\lambda|^{-N_0}$ on $S$. Thus $(5.4)$
is bounded by
$$C|\lambda|^{{N_0\over 2}Re(z)}\int_S |\bar{H}(x)^{z \over 4} - \bar{H}_N(x)^{z \over 4}|\eqno (5.5)$$
Let $N_1$ be such that the measure of $\{x: |\bar{H}(x)| < |\lambda|^{-N_1}\}$ is less than $C|\lambda|
^{-1 - {N_0\over 2}Re(z)}$. Then removing this set from the domain of integration of $(5.5)$ will change
$(5.5)$ by at most $C|\lambda|^{-1}$. Since this is smaller than the right-hand side of $(5.1a)$, it suffices to find upper bounds for
$$C|\lambda|^{{N_0\over 2}Re(z)}\int_{\{x \in S: |\bar{H}(x)| \geq |\lambda|^{-N_1}\}}|\bar{H}(x)^{z \over 4} - 
\bar{H}_N(x)^{z \over 4}| \eqno (5.5')$$
By taking $N$ to be sufficiently large and the constant we called $c$ in the definition of the subrectangles
to be sufficiently small, we can ensure that $|\bar{H}(x)^{z \over 4} - \bar{H}_N(x)^{z \over 4}|$
is bounded by $C'|\lambda|^{-1 - {N_0\over 2}Re(z)}$ whenever $|\bar{H}(x)| > |\lambda|^{-N_1}$. Hence in this case
$(5.5')$ is bounded by $C''|\lambda|^{-1}$, better than what we need. Hence we may replace $H^{**}(x)$ by
$H^{**}_N(x)$ in our subsequent arguments. It is worth noting that here the constants $C'$ and $C''$ do
depend linearly on $|Im(z)|$. 

One polynomializes the $D(x)$ factor in much the same way as we dealt with $H^{**}(x)$. The conclusion
is that if $D_N(x)$ denotes the determinant of $H_N(x)$ in the $u$ and $v$ variables, then
if $N$ is chosen sufficiently large and $c$ is chosen sufficiently small, then one can replace $D(x)$ and
$D_N(x)$ in our arguments without affecting the conclusions. 

In summary, to show $|G_z^S(\lambda)|$ is bounded by the right-hand side of $(5.1a)$, thereby proving 
Theorem 5.1, we must show that if $N$ is sufficiently
large then  $|(G_z^S)_N(\lambda)|$ satisfies the same bounds, where 
$$(G_z^S)_N(\lambda) = \int e^{-i\lambda_1 x_1... - i\lambda_n x_n  - i\lambda_{n+1} H_N(x)}|D_N(x)|^{\delta z}$$
$$|H^*(x)|^{-{z \over \max(d,2)}} |H^{**}_N(x)|^z \Psi(x)\,dx_1...\,dx_n \eqno (5.6)$$
Analogous to in section 3, if the maximal $|\lambda_i|$ is not $|\lambda_{n+1}|$ then by integrating by parts 
in the $x_i$ variable one can obtain stronger decay than is needed. So without loss of 
generality we henceforth assume $|\lambda_{n+1}| \geq |\lambda_i|$ for all $i < n+1$. 

\noindent We next show that there is a fixed $C$ depending only on the function $H(x)$ such that 
$$|H^*(x)|^{-{1 \over \max(d,2)}} |H^{**}_N(x)| < C \eqno (5.7)$$
In particular, the damping factors are uniformly bounded for a fixed value of $Re(z)$.
To see why $(5.7)$ holds,
recall that by Lemma 4.3, $N(\bar{H}) = 2N(H)$. Hence by the corollary to Theorem 3.2, there is a $C_1$ 
depending only on $H$ such that $|\bar{H}^{**}(x)| < C_0|H^*(x)|^2$ and therefore $|H^{**}(x)| < 
C_1|H^*(x)|^{1 \over 2}$. (Technically, this corollary was 
proven for the real-analytic case only, but the smooth cases under consideration can be covered by considering
a truncation of the Taylor series to high enough order and observing that the error term is easily bounded
in such cases by a constant times $|H^*(x)|^{1 \over 2}$). Hence we have
$$|H^*(x)|^{-{1 \over 2}} |H^{**}(x)| < C_2\eqno (5.8a)$$
Recall we are assuming that $S$ comes from an $R$ for which $\inf_R |H^*(x)| > |\lambda|^{-N_0}$ for some $N_0$. 
Thus if $N$ were chosen sufficiently large and $c$ sufficiently small so that $|\bar{H}(x) -  \bar{H}_N(x)|
< C_2^4|\lambda|^{-2N_0}$, then $\big||H^*(x)|^{-2}|\bar{H}(x)| -  |H^*(x)|^{-2}|\bar{H}_N(x)|)\big| < 
C_2^4$ and thus in view of $(5.8a)$ we have 
$$|H^*(x)|^{-2} |H^{**}_N(x)|^4 < 2C_2^4 $$
In turn, this implies
$$|H^*(x)|^{-{1 \over 2}} |H^{**}_N(x)| < 2^{1 \over 4}C_2 \eqno (5.8b)$$
Since $|H^*(x)|^{-{1 \over \max(2,d)}} \leq |H^*(x)|^{-{1 \over 2}}$, this establishes $(5.7)$ as needed.

We now proceed to show that $(G_z^S)_N(\lambda)$ is bounded by the right hand side of 
$(5.1a)$. We split the integral into $(G_z^S)_N(\lambda) = (H_z^S)_N(\lambda) + (I_z^S)_N(\lambda)$, where $(H_z^S)_N
(\lambda)$ is the part of $(5.6)$ over where $|D_N(x)| > |\lambda|^{-\delta_2}$ and $(I_z^S)_N(\lambda)$
is the integral over where $|D_N(x)| < |\lambda|^{-\delta_2}$. Here $\delta_2$ is a constant to be 
determined by our arguments, unrelated to $\delta_1$ or $\delta$, which will only depend on the dimension 
$n$. 

We start with the analysis of $(H_z^S)_N(\lambda)$. On the domain of $(H_z^S)_N(\lambda)$, the determinant of 
$H_N$ in the $u$ and $v$ variables is bounded below by something relatively large, and methods
similar to those used to deal with 2-d oscillatory integrals with nonvanishing Hessian will be used in our
analysis. We divide first divide $S$ into identical cubes $C_i$ of radius $c'|\lambda|^{-{\delta_2}}$, where
$c'$ will depend on $N$ and $H(x)$. We examine the contribution to $(H_z^S)_N(\lambda)$ from each cube $C_i$. The contribution
is nonzero only if $|D_N(x)| > |\lambda|^{-\delta_2}$ for some $x$ in $C_i$, so we assume this in fact the
case. Assuming $c'$ were chosen sufficiently small, on $C_i$ we have that
$$|D_N(x)| > {1 \over 2}|\lambda|^{-\delta_2} $$
Recalling that $u$ and $v$ are the directions in which the Hessian determinant $D_N(x)$ are taken, the $u$
derivative of the phase function $-\lambda_1 x_1-... - \lambda_n x_n  - \lambda_{n+1} H_N(x)$ of $(H_z^S)_N(\lambda)$
can be rewritten as $-\lambda_{n+1}(\partial_u H_N(x)  - a_1)$ for some constant $a_1$. One can similarly 
rewrite the $v$ derivative as $-\lambda_{n+1}(\partial_v H_N(x)  - a_2)$.
We write the contribution to $(H_z^S)_N(\lambda)$ coming from $C_i$ as $J_1 + J_2 + J_3$, where 
$$J_1 = \int_{\{x:|\partial_u H_N(x) - a_1| > |\lambda|^{-{1 \over 3}}\}} e^{-i\lambda_1 x_1-... - i\lambda_n x_n  - i\lambda_{n+1} H_N(x)}$$
$$\times |D_N(x)|^{\delta z} |H^*(x)|^{- {z \over \max(2,d)}} |H^{**}_N(x)|^z\Psi(x)\,dx_1...\,dx_n \eqno (5.9a)$$ 
$$J_2 = \int_{\{x: |\partial_u H_N(x)  - a_1| < |\lambda|^{-{1 \over 3}},\,\,|\partial_v H_N(x) - a_2|> |\lambda|^{-{1 \over 3}}\}} 
e^{-i\lambda_1 x_1-... - i\lambda_n x_n  - i\lambda_{n+1} H_N(x)}$$
$$\times |D_N(x)|^{\delta z} |H^*(x)|^{- {z \over \max(2,d)}}|H^{**}_N(x)|^z\Psi(x)\,dx_1...\,dx_n \eqno (5.9b)$$ 
$$J_3 = \int_{\{x: |\partial_u H_N(x) - a_1| < |\lambda|^{-{1 \over 3}},\,\,|\partial_v H_N(x) - a_2| < |\lambda|^{-{1 \over 3}}\}} 
e^{-i\lambda_1 x_1-... - i\lambda_n x_n  - i\lambda_{n+1} H_N(x)}$$
$$\times |D_N(x)|^{\delta z} |H^*(x)|^{-{z \over \max(2,d)}}|H^{**}_N(x)|^z\Psi(x)\,dx_1...\,dx_n \eqno (5.9c)$$ 
We start with $J_1$. We integrate by parts in the $u$ direction integrating $\lambda_{n+1}(\partial_u H_N(x)  - a_1)
e^{-i\lambda_1 x_1... - i\lambda_n x_n  - i\lambda_{n+1} H_N(x)}$ and differentiating ${1 \over \lambda_{n+1}
(\partial_u H_N(x)  - a_1)}$ times the rest of the integrand. The derivative can land on ${1 \over \lambda_{n+1}
(\partial_u H_N(x)  - a_1)}$ or any of the various factors appearing in $(5.9a)$. In all cases, one takes 
absolute values and integrates as in the proof of the Van der Corput lemma. The fact that all the factors
are polynomials, as is $\partial_u H_N(x)  - a_1$, ensures that one integrates over boundedly many intervals
on each of which the relevant derivative is monotone. Hence one can integrate back each derivative, just like
in the proof of the Van der Corput lemma. As a result, $(5.9a)$ is bounded by a constant times the supremum of 
$|{1 \over \lambda_{n+1}(\partial_u H_N(x)  - a_1)}|$ on its domain, which is at most $C|\lambda|^{-{2 \over 3}}$. 
It is worth
mentioning that $C'$ here actually depends linearly on $|Im(z)|$ due to the terms where the derivative lands on
the damping factors. Although the integrand itself is uniformly bounded for fixed $Re(z)$, the derivatives
landing on damping factors incur factors linear in $z$ and the process of taking absolute values gives the 
linear dependence on $|Im(z)|$. 

One bounds $J_2$ exactly in the same way as $J_1$, using the $v$ derivative in place of the $u$ derivative. 
As for $J_3$, the fact that $|D_N(x)| > {1 \over 2}|\lambda|^{-\delta_2}$ on $C_i$ implies that the gradients of
$\partial_u H_N$ and $\partial_v H_N$ are both of magnitude at least $C|\lambda|^{-\delta_2}$ on $C_i$. Hence 
if the constant $c'$ in the definition $c'\delta_2$ of the radius of $C_i$ is small enough, the level sets of
$\partial_u H_N$ and $\partial_v H_N$ are smooth manifolds (don't "self-intersect"). In particular if $y_3,...,y_n$ 
denote orthonormal
directions perpendicular to the plane generated by $u$ and $v$, the coordinate change from $(u,v,y_3,...,y_n)$
to $(\partial_u H_N,\partial_v H_N,y_3,...,y_n)$ is well-defined and has Jacobian bounded below by ${1 \over 2}
|\lambda|^{-\delta_2}$. In particular, the inverse image of the points where $|\partial_u H_N  - a_1|, |\partial_v H_N - a_2| 
< |\lambda|^{-{1 \over 3}}$ under this map has measure at most $|\lambda|^{-{2 \over 3}} \times
2|\lambda|^{\delta_2} \leq 2|\lambda|^{\delta_2 - {2 \over 3}}$. Hence $J_3$ is bounded by the supremum of
the integrand times this, or $C|\lambda|^{\delta_2 - {2 \over 3}}$. 

Adding up $J_1 + J_2 + J_3$, we see that the contribution to $(H_z^S)_N(\lambda)$ from the cube $C_i$ is at most
$C'|\lambda|^{\delta_2 - {2 \over 3}}$. There are at most $C''|\lambda|^{n\delta_2}$ such cubes, so we conclude that
$$|(H_z^S)_N(\lambda)| \leq C|\lambda|^{(n+1)\delta_2 - {2 \over 3}} \eqno (5.10a)$$
Thus choosing $\delta_2$ small enough depending on $n$ only gives 
$$|(H_z^S)_N(\lambda)| \leq C|\lambda|^{3 \over 5} \eqno (5.10b)$$
This is better than the estimate we need. We now proceed to bounding $(I_z^S)_N$, given by
$$(I_z^S)_N = \int_{\{x: |D_N(x)| < |\lambda|^{-\delta_2}\}} e^{-i\lambda_1 x_1- ... - i\lambda_n x_n  - 
i\lambda_{n+1} H_N(x)}|D_N(x)|^{\delta z}$$
$$\times  |H^*(x)|^{-{z \over \max(d,2)}} |H^{**}_N(x)|^z \Psi(x)\,dx_1...\,dx_n \eqno (5.11)$$
Recall by $(5.8b)$ there is some constant $C_0$ depending only $H(x)$ such that $|H^{**}_N(x)| < 
C_0|H^*(x)|^{{1 \over 2}}$. Thus we may write $(I_z^S)_N = \sum_{i = 0}^{\infty} A_i$, where $A_i$ is
the portion of $(5.11)$ where $2^{-i - 1}C_0|H^*(x)|^{{1 \over 2}} \leq |H^{**}_N(x)| < 
2^{-i}C_0|H^*(x)|^{{1 \over 2}}$. Next, recall that $S$ was defined as a subrectangle of some dyadic rectangle
$\prod_{j = 1}^n [2^{-k_j - 1}, 2^{-k_j}]$. We do a change of variable in the integral $(5.11)$, turning
$x$ into $2^{-k}x = (2^{-k_1}x_1,...,2^{-k_n}x_n)$. Then $A_i$ can be written as 
$$A_i = 2^{-\sum_j k_j}\int_{\{x \in B_i: |D_N(2^{-k}x)| < |\lambda|^{-\delta_2}\}} e^{-i2^{-k_1}\lambda_1 
x_1-... - i 2^{-k_n}\lambda_nx_n  - i\lambda_{n+1} H_N(2^{-k}x)}$$
$$\times |D_N(2^{-k}x)|^{\delta z} |H^*(2^{-k}x)|^{-{z \over \max(d,2)}} |H^{**}_N(2^{-k}x)|^z \Psi(2^{-k}x)
\,dx_1...\,dx_n \eqno (5.12)$$
Here $B_i$ are the $x$ in $S$ such that $2^{-i - 1}C_0|H^*(2^{-k}x)|^{{1 \over 2}} \leq 
|H^{**}_N(2^{-k}x)| < 2^{-i}C_0|H^*(2^{-k})|^{{1 \over 2}}$. 

Next, note that there exists a constant $C_1$ such that
$(H^{**}_N(2^{-k}x))^2$ is within a factor of $C_1$ of $\sum_{i,j}|\partial_{x_ix_j}^2\big(H_N(2^{-k}x)\big)|$,
which in turn is within a constant factor of $\sup_{i,j}|\partial_{x_ix_j}^2\big(H_N(2^{-k}x)\big)|$. 
It is well known (see [St2] p. 343) that there are a finite set of directions $\xi_1,...,\xi_p$ such that
every second partial operator $\partial_{x_ix_j}^2$ can be written as a linear combination of the
$(\xi_l \cdot \nabla)^2$. In particular, this means that $(H^{**}_N(2^{-k}x))^2$ is within a constant factor
of $\sup_l |(\xi_l \cdot \nabla)^2 \big(H_N(2^{-k}x)\big)|$. 

We now write $(5.12)$ as $\sum_{l = 1}^p D_{il}$, where 
$D_{il}$ denotes the portion of $A_i$ where $|(\xi_l \cdot \nabla)^2\big(H_N(2^{-k}x)\big)|$ is greater than $|(\xi_m
\cdot \nabla)^2 \big(H_N(2^{-k}x)\big)|$ for $m \neq l$. (In the unlikely event two or more of these functions are 
the same we do not repeat any functions). Thus throughout the domain of integration of $D_{il}$, one has that 
$|(\xi_l \cdot \nabla)^2\big(H_N(2^{-k}x)\big)|$ is within a constant factor of $(H^{**}_N(2^{-k}x))^2$, which in 
turn within a constant factor of $2^{-2i}|H^*(2^{-k}x)|$. Call this domain
of integration $E_{il}$. Then
$$D_{il} = 2^{-\sum_j k_j}\int_{E_{il}}e^{-i2^{-k_1}\lambda_1 
x_1-... - i 2^{-k_n}\lambda_nx_n  - i\lambda_{n+1} H_N(2^{-k}x)}|D_N(2^{-k}x)|^{\delta z}$$
$$\times  |H^*(2^{-k}x)|^{-{z \over \max(d,2)}} |H^{**}_N(2^{-k}x)|^z \Psi(2^{-k}x)
\,dx_1...\,dx_n \eqno (5.13)$$
Let $P(x)$ denote the phase function in $(5.13)$. Then since $P(x)$ differs from $\lambda_{n+1}H_N(2^{-k}x)$
by a linear function, there are constants $C$ and $C'$, depending only on $H(x)$, such that on $E_{il}$ we have
$$C 2^{-2i}|\lambda_{n+1}H^*(2^{-k}x)| \leq |(\xi_l \cdot \nabla)^2 P(x)| \leq C'2^{-2i}|\lambda_{n+1}
H^*(2^{-k}x)| \eqno (5.14)$$
Recall that $H^*(2^{-k}x)$ is within a constant factor of a fixed value $H^*(R)$ on the 
dyadic rectangle $R$ which $S$ is a part of. Hence instead of $(5.14)$ we can use
$$C 2^{-2i}|\lambda_{n+1}H^*(R)| \leq |(\xi_l \cdot \nabla)^2 P(x)| \leq C'2^{-2i}|\lambda_{n+1}
H^*(R)| \eqno (5.14')$$
Furthermore, since $\lambda_{n+1}$ is assumed to be larger than $\lambda_i$ for $i \neq n+1$, we can also use
$$C''2^{-2i}|\lambda| |H^*(R)| \leq |(\xi_l \cdot \nabla)^2 P(x)| \leq C'''2^{-2i}|\lambda| |H^*(R)| 
\eqno (5.14'')$$
In view of $(5.14'')$, we will now argue as in the proof of the Van der Corput lemma in the $\xi_l$ 
direction.
We break $(5.13)$ up as $D_{il}^1 + D_{il}^2$, where $D_{il}^1$ is the portion where $|(\xi_l \cdot \nabla) P(x)|$
is less than $2^{-i}|\lambda|^{{1 \over 2}}|H^*(R)|^{{1 \over 2}}$ and $D_{il}^2$ is the portion where $|(\xi_l \cdot
\nabla) P(x)|$ is greater than $2^{-i}|\lambda|^{{1 \over 2}}|H^*(R)|^{{1 \over 2}}$. 

To estimate $D_{il}^1$, we take absolute values and integrate. In view of $(5.14'')$, $D_{il}^1$ is bounded
by $2^i|\lambda|^{-{1 \over 2}} |H^*(R)|^{-{1 \over 2}}$ times the supremum of the absolute value of the integrand.
In the integrand, $\big||D_N(2^{-k}x)|^{\delta z}\big| = |D_N(2^{-k}x)|^{\delta Re(z)}$, which in view of 
the definition of the
domain of integration of $(5.12)$ is at most $|\lambda|^{-\delta \delta_2 Re(z)}$. Next, the factor  
$|H^{**}_N(2^{-k}x)|^z$ has magnitude $|H^{**}_N(2^{-k}x)|^{Re(z)} \leq C2^{-iRe(z)}|H^*(2^{-k}x)|^
{{1 \over 2}Re(z)} \leq  C'2^{-iRe(z)}|H^*(R)|^{{1 \over 2}Re(z)}$. Thus $|D_{il}^1|$ is at most a constant
times
$$2^i|\lambda|^{-{1 \over 2}}|H^*(R)|^{-{1 \over 2}}\big(2^{-\sum_j k_j} \times |\lambda|^{-\delta \delta_2 
Re(z)}\times |H^*(R)|^{-{1 \over \max(2,d)} Re(z)}\times 2^{-iRe(z)}|H^*(R)|^{{1 \over 2}Re(z)}\big) $$
$$ = C2^{-\sum_j k_j} |\lambda|^{-{1 \over 2} - \delta \delta_2 Re(z)}2^{-i(Re(z) - 1)}|H^*(R)|^
{Re(z)({1 \over 2} - {1 \over \max(d,2)})-{1 \over 2}}\eqno (5.15)$$
These will be seen to be the bounds we need. We now move on to $D_{il}^2$. For this, we integrate by parts in the
$\xi_l$ direction as in the proof of the Van der Corput lemma, writing $e^{iP(x)} = ((\xi_l \cdot \nabla) P(x)
e^{iP(x)}) / ((\xi_l \cdot \nabla) P(x))$ and integrating the numerator. The derivative can land on several
factors. For each term thus generated, we take absolute values and integrate. The polynomial character of
all the factors (except for $\Psi(2^{-k}x)$, which doesn't cause any problems) as well as $P(x)$ and the 
functions
defining the domain of integration ensures that we integrate over boundedly many intervals on each of
which the differentiated factor is monotone. Hence the proof of the Van der Corput applies and $|D_{il}^2|$
is at most the supremum of $|(\xi_l \cdot \nabla) P(x)|^{-1}$ on the domain of integration times the
supremum of the magnitude of the integrand on the domain of the integration. The same bounds hold for 
the endpoint terms coming from the integration by parts. 

Note that on the domain of integration of $D_{il}^2$, $|(\xi_l \cdot \nabla) P(x)|^{-1} < 
2^i|\lambda H^*(R)|^{-{1 \over 2}}$. In the analysis of $D_{il}^1$, we bounded $|D_{il}^1|$ by this same
factor times the supremum of the magnitude of the integrand. Hence like before we have
$$|D_{il}^2| \leq C2^{-\sum_j k_j} |\lambda|^{-{1 \over 2} - \delta \delta_2 Re(z)}2^{-i(Re(z) - 1)}|H^*(R)|^
{Re(z)({1 \over 2} - {1 \over \max(d,2)})-{1 \over 2}}\eqno (5.16)$$
Adding $(5.15)$ to $(5.16)$, and then summing over all $l$ gives the following, where $A_i$ was given by 
$(5.12)$.
$$|A_i| \leq C'2^{-\sum_j k_j} |\lambda|^{-{1 \over 2} - \delta \delta_2 Re(z)}2^{-i(Re(z) - 1)}|H^*(R)|^
{Re(z)({1 \over 2} - {1 \over \max(d,2)})-{1 \over 2}}\eqno (5.17)$$
Since we are assuming $Re(z) > 1$, if we add $(5.17)$ over all $i$ we get
$$|(I_z^S)_N| \leq C''2^{-\sum_j k_j} |\lambda|^{-{1 \over 2} - \delta \delta_2 Re(z)}|H^*(R)|^
{Re(z)({1 \over 2} - {1 \over \max(d,2)})-{1 \over 2}}\eqno (5.18)$$
Next, note that $2^{-\sum_j k_j}$ is the area of the rectangle $R$ that $S$ is a part of, and that $H^*(x)$
is within a constant factor of $H^*(R)$ on $R$. Hence $(5.18)$ can be reexpressed as 
$$|(I_z^S)_N| \leq C'''|\lambda|^{-{1 \over 2} - \delta \delta_2 Re(z)} \int_R |H^*(x)|^
{Re(z)({1 \over 2} - {1 \over \max(d,2)})-{1 \over 2}} \eqno (5.19)$$
As mentioned before, by Theorem 1.2 of [G1], the supremum of the $t$ for which $|H^*(t)|^{-t}$ is integrable
on a neighborhood of the origin is ${1 \over d}$. Since $Re(z) > 1$, the exponent in $(5.19)$ is greater
than $-{1 \over \max(d,2)}$ which is itself at least $-{1 \over d}$. Hence the function in $(5.19)$ is 
integrable and we have
$$|(I_z^S)_N| \leq C'''|\lambda|^{-{1 \over 2} - \delta \delta_2 Re(z)} \eqno (5.20)$$
Adding this to our estimate $(5.10b)$ for $|(H_z^S)_N|$, we see that $(G_z^S)_N$ defined in $(5.6)$ 
satisfies
$$|(G_z^S)_N| \leq C_4|\lambda|^{-{1 \over 2} - \delta \delta_2 Re(z)} \eqno (5.21)$$
Recall that in order to prove Theorem 5.1, we had to show that $|(G_z^S)_N|$ was bounded by the right-hand
side of $(5.1a)$ for sufficiently small $\delta$. Since $\delta_2$ was a constant depending only on
the dimension $n$, and all constants $C, C'$ etc appearing here grow at most linearly in $|Im(z)|$, $(5.21)$
completes the proof of Theorem 5.1.

\noindent {\bf The proofs of Theorems 1.1 and 1.3.}

Using the version of analytic interpolation for maximal operators on the $e^{z^2}\sigma_z$ (see p.482 of
[St2] for another example of this), we now use Lemma 4.4 and Theorem 5.1 to complete the proofs of Theorem 1.1 and
1.3. Suppose $M_z$ is uniformly bounded on $L^{\infty}$ on any $Re(z) = c$ with $c > -a$, and is uniformly bounded on 
$L^2$ on any $Re(z) = c$ with $c > 1$. Then $M_0 = M$ is bounded on $L^p$ for $p > p_0$,
where ${1 \over p_0} = {1 \over a + 1}{1 \over \infty} + {a \over a + 1}{1 \over 2} = {a \over 2a + 2}$.
Hence $p_0 = {2a + 2 \over a}$. 

We first suppose we are in the setting of Theorem 1.1a) and the corresponding part of Theorem 1.3. Then 
$d \leq 2$, and the number $M = \max(2,m)$ of Lemma 4.4 is 2. Hence we are in case c) of Lemma 4.4 and
$a$ can be taken as large as one wants. Thus we
concude that $p_0$ can be taken to be $\lim_{a \rightarrow \infty} {2a + 2 \over a} = 2$ as required.

We next move to the setting of Theorem 1.1b) and the corresponding part of Theorem 1.3. Now $d > 2$ and
$M = \max(2,m) \leq d$. Hence by case a) of Lemma 4.4, we can take $a = {2 \over d - 2}$. Hence $p_0 = {{4 \over d - 2} + 2 
\over {2 \over d - 2}} = d$ as required. The sharpness of $p_0$ when $y \notin T_y(S)$ holds because
by Theorem 1.2b) of [G1], $\int |H|^{-{1 \over d}}$ is infinite on any neighborhood of the origin and
by [IoSa1] this is sufficient for unboundedness of $M$ on $L^d$ so long as $y \notin T_y(S)$.

Lastly, suppose we are in the setting of Theorem 1.1c) or the
corresponding part of Theorem 1.3. Then $M =  m > \max(2,d)$. 
If $d \leq 2$, then since $M > 2$ case b) of Lemma 4.4 says $a = {2 \over M - 2}$. In 
this case $p_0 = {{4 \over M - 2} + 2 \over {2 \over M - 2}} = M$ as needed. If $d > 2$, then since $M > d$ 
case a) of Lemma 4.4 says that $a = {2 \over M - 2}$ and thus once again $p_0 = M$ as needed. This completes
the proofs of Theorems 1.1 and 1.3.

\noindent {\bf 6. References.}

\noindent [Bo] J. Bourgain, {\it Averages in the plane over convex curves and maximal operators},
J. Anal. Math. {\bf 47} (1986), 69--85. \parskip = 4pt\baselineskip = 4pt

\noindent [BrNaWa] J. Bruna, A. Nagel, and S. Wainger, {\it Convex hypersurfaces and Fourier transforms},
Ann. of Math. (2) {\bf 127} no. 2, (1988), 333--365. 

\noindent [Ch] M. Christ, {\it Hilbert transforms along curves. I. Nilpotent groups}, Annals of Mathematics
(2) {\bf 122} (1985), no.3, 575-596.

\noindent [CoDiMaMu] M. Cowling, S. Disney, G Mauceri, and D. Muller {\it Damping oscillatory integrals}, 
Invent. Math. {\bf 101}  (1990),  no. 2, 237--260.

\noindent [CoMa] M. Cowling, G. Mauceri, {\it Oscillatory integrals and Fourier transforms of surface 
carried measures},  Trans. Amer. Math. Soc. {\bf 304} (1987),  no. 1, 53--68. 

\noindent [G1] M. Greenblatt, {\it Oscillatory integral decay, sublevel set growth, and the Newton
polyhedron}, Math. Annalen {\bf 346} (2010), no. 4, 857-895.

\noindent [G2] M. Greenblatt, {\it A Coordinate-dependent local resolution of singularities and 
applications},  J. Funct. Anal.  {\bf 255}  (2008), no. 8, 1957-1994.

\noindent [Gr] A. Greenleaf, {\it Principal curvature and harmonic analysis}, 
Indiana Univ. Math. J. {\bf 30} (1981), no. 4, 519--537. 

\noindent [IkKeMu1] I. Ikromov, M. Kempe, and D. M\"uller, {\it Damped oscillatory integrals and boundedness of
maximal operators associated to mixed homogeneous hypersurfaces} (English summary) Duke Math. J. {\bf 126} 
(2005), no. 3, 471--490.

\noindent [IkKeMu2] I. Ikromov, M. Kempe, and D. M\"uller, {\it Estimates for maximal functions associated
to hypersurfaces in $R^3$ and related problems of harmonic analysis}, to appear, Acta Math.

\noindent [Io]  A. Iosevich, {\it Maximal operators associated to families of flat curves in the plane},
Duke Math. J. {\bf 76} no. 2 (1994) 633-644. 

\noindent [IoSa1] A. Iosevich, E. Sawyer, {\it Oscillatory integrals and maximal averages over homogeneous
surfaces}, Duke Math. J. {\bf 82} no. 1 (1996), 103-141.

\noindent [IoSa2] A. Iosevich, E. Sawyer, {\it Maximal averages over surfaces},  Adv. Math. {\bf 132} 
(1997), no. 1, 46--119.

\noindent [IoSaSe] A. Iosevich, E. Sawyer, A. Seeger, {\it On averaging operators associated with convex 
hypersurfaces of finite type}, J. Anal. Math. {\bf 79} (1999), 159--187

\noindent [MocSeSo] G. Mockenhaupt, A. Seeger, and C. Sogge, {\it Wave front sets, local smoothing and 
Bourgain's circular maximal theorem}, Ann. of Math. (2) {\bf 136} (1992), no. 1, 207--218.

\noindent [NaSeWa] A. Nagel, A. Seeger, and S. Wainger, {\it Averages over convex hypersurfaces},
Amer. J. Math. {\bf 115} (1993), no. 4, 903--927.

\noindent [PhSt] D. H. Phong, E. M. Stein, {\it The Newton polyhedron and
oscillatory integral operators}, Acta Mathematica {\bf 179} (1997), 107-152.

\noindent [Sc] H. Schulz, {\it Convex hypersurfaces of finite type and the asymptotics of their Fourier 
transforms}, Indiana Univ. Math. J. {\bf 40} (1991), no. 4, 1267--1275. 

\noindent [So] C. Sogge, {\it Maximal operators associated to hypersurfaces with one nonvanishing principal 
curvature}  (English summary)  in {\it Fourier analysis and partial differential equations} (Miraflores de 
la Sierra, 1992),  317--323, Stud. Adv. Math., CRC, Boca Raton, FL, 1995. 

\noindent [SoSt] C. Sogge and E. Stein, {\it Averages of functions over hypersurfaces in $R^n$}, Invent.
Math. {\bf 82} (1985), no. 3, 543--556.

\noindent [St1] E. Stein, {\it Maximal functions. I. Spherical means.} Proc. Nat. Acad. Sci. U.S.A. 
{\bf 73} (1976), no. 7, 2174--2175. 

\noindent [St2] E. Stein, {\it Harmonic analysis; real-variable methods,
orthogonality, and oscillatory integrals}, Princeton Mathematics Series Vol. 
43, Princeton University Press, Princeton, NJ, 1993.

\noindent [V] A. N. Varchenko, {\it Newton polyhedra and estimates of oscillatory integrals}, Functional 
Anal. Appl. {\bf 18} (1976), no. 3, 175-196.

\line{}
\line{}

\noindent Department of Mathematics, Statistics, and Computer Science \hfill \break
\noindent University of Illinois at Chicago \hfill \break
\noindent 322 Science and Engineering Offices \hfill \break
\noindent 851 S. Morgan Street \hfill \break
\noindent Chicago, IL 60607-7045 \hfill \break
\end